\documentclass[11pt]{article}
\usepackage[margin=1in]{geometry}                
\usepackage[parfill]{parskip}   
\usepackage{graphicx}
\usepackage{amssymb}
\usepackage{epstopdf}
\usepackage{amsmath}
\usepackage{amsthm}
\usepackage{fancyhdr}

\newtheorem{theorem}{Theorem}[section]
\newtheorem{lemma}[theorem]{Lemma}
\newtheorem{proposition}[theorem]{Proposition}
\newtheorem{remark}[theorem]{Remark}
\numberwithin{equation}{section}

\begin{document}

\title{Hydrodynamic limit of the zero range process on a randomly oriented graph}
\author{M\'arton Bal\'azs\thanks{University of Bristol; \texttt{m.balazs@bristol.ac.uk}}
 \and
Felix Maxey-Hawkins\thanks{University of Bristol; \texttt{fm17392@bristol.ac.uk}}}


\maketitle

\begin{abstract}
We prove the hydrodynamic limit of a totally asymmetric zero range process on a torus with two lanes and randomly oriented edges. The asymmetry implies that the model is non-reversible. The random orientation of the edges is constructed in a bistochastic fashion which keeps the usual product distribution stationary for the quenched zero range model. It is also arranged to have no overall drift along the $\mathbb Z$ direction, which suggests diffusive scaling despite the asymmetry present in the dynamics. Indeed, using the relative entropy method, we prove the quenched hydrodynamic limit to be the heat equation with a diffusion coefficient depending on ergodic properties of the orientation of the edges.

 The zero range process on this graph turns out to be non-gradient. Our main novelty is the introduction of a local equilibrium measure which decomposes the vertices of the graph into components of constant density. A clever choice of these components eliminates the non-gradient problems that normally arise during the hydrodynamic limit procedure.
\end{abstract}

\section{Introduction}

Hydrodynamic limits of interacting particle systems have a long and rich history, starting from relatively simpler cases where the model is of \emph{gradient type} and the dynamics is reversible for the stationary distribution, to more complicated setups with non-gradient models and/or non-reversible dynamics. Models can be mean zero or have drift, which generally decides between Eulerian scaling and hyperbolic limiting PDE, or diffusive scaling with parabolic limit. We refer the reader to the fundamental book of Kipnis and Landim \cite{landim} for overview and details.

Extra complications arise when the dynamics is run in a random environment. We consider such a scenario by running a totally asymmetric nearest neighbour zero range process on a graph with vertices in $\mathbb Z\times\{-1,1\}$ and randomly oriented edges (i.e., \emph{environment}). We call the process \emph{totally asymmetric} in a local sense as particles cannot move against the arrows. As detailed shortly below, the graph is \emph{bistochastic}, which keeps the usual product distribution stationary for zero range. It also has zero drift in the $\mathbb Z$ direction, which comes with diffusive scaling. Together with some ergodicity assumptions on the environment, this puts us in the situation that our model
\begin{itemize}
 \item lives in a random environment,
 \item has an i.i.d.\ stationary distribution,
 \item turns out to be non-gradient,
 \item is non-reversible, in fact totally asymmetric,
 \item has, nevertheless, diffusive scaling and the heat equation as its (quenched) hydrodynamics.
\end{itemize}
This last statement is what we prove in this paper. The heat equation comes with a coefficient that only depends on ergodic properties of the environment. Of the several routes in the literature, we follow the relative entropy method worked out by Yau \cite{yau_relentr}. Not having the gradient property generally poses significant difficulties in the method.

The main novelty of this paper is proving hydrodynamics in a non-gradient, non-reversible random environment. This is  made possible by a clever choice of the local equilibrium measure that eliminates the difficulties coming from the lack of the gradient property. This is done by carefully looking at the orientation of the underlying graph. The observation paves the way for conventional hydrodynamic arguments, hence providing a reasonably simple proof in the above setup that usually requires elaborate arguments like sector conditions or the two block estimate. We believe our method could be extended to a randomly oriented graph with more than two lanes subject to some restrictions on the orientation. This might provide some generalisations but a full treatment of all divergence free environments seems out of reach by our methods. 

The usual non-gradient method introduced by Quastel \cite{q2} and Varadhan \cite{v2} and also treated in Kipnis and Landim \cite{landim} is an adaptation of the entropy method. 
However, in this paper we follow the previous work on non-gradient models by Funaki, Uchiyama and Yau \cite{funaki} for a reversible model and later Komoriya \cite{komoriya} whose model was non-reversible. These authors adapted the relative entropy method to non-gradient models by introducing the local equilibrium state of second order approximation. These local equilibrium states were constructed differently from ours because the jump rates were non-random and translation invariant, but they are the inspiration for our construction.

We mention some further results from the literature concerning random environments. Random walks in random environments have been studied extensively and a great variety of results are available. Bistochastic, non-reversible environments such as that of Kozma and T\'oth \cite{kozma} are an important development in recent years, generalising more classical, reversible scenarios. Our work is motivated by the desire to study interacting particle systems in analogous random environments.

 Koukkous \cite{koukkous} used the entropy method to prove that the heat equation is the hydrodynamic limit of a symmetric zero range process on a $d$-dimensional torus for which the sites have random jump rates and the model is therefore non-gradient.
 Faggionato and Martinelli \cite{faggionato} used the non-gradient method to prove that the hydrodynamic limit of an exclusion process in dimension $d \geq 3$ with random transition rates and satisfying a detailed balance condition is the heat equation. Quastel \cite{quastel} proved the same result in all dimensions using a variation of the non-gradient method. 
 Goncalves and Jara \cite{gj} obtained the hydrodynamic limit of a zero range process on the $d$-dimensional torus with the jump rates across edges given by a random environment, making use of the entropy method combined with a corrected empirical measure and homogenisation results. Faggionato \cite{f2} extended the result to a supercritical percolation cluster.
 We notice that the dynamics in these models are reversible, as opposed to our setup. Bahadoran et al. \cite{bahadoran} provides a review of work on one dimensional asymmetric zero range processes in random environments, including hydrodynamic limits. 

Next we describe the randomly oriented graph our zero range process lives on. It is one of the simplest graphs that can be given a non-trivial bistochastic orientation. The hydrodynamic limit is performed on the torus $\mathbb{T}_N \times \{-1,1\}$, where $\mathbb{T}_N = \mathbb{Z}/N\mathbb{Z}$. We interpret the first part of the Cartesian product as horizontal coordinates and the second part as distinguishing between \emph{lower} and \emph{upper} vertices. We will call each pair of vertices $(j,1),(j,-1)\in\mathbb{T}_N\times\{-1,1\}$ a \emph{site}. Each vertex is connected by horizontal and diagonal edges to both vertices on its left neighbouring site as well as to both vertices on its right neighbouring site.

The edges are oriented subject to these rules:
\begin{itemize}
 \item The graph is bistochastic, in other words divergence-free. That translates to our case as each vertex having exactly two in-edges and two out-edges. Considering the directed graph as a flow this would exactly translate to this flow being divergence free everywhere.
 \item There is no overall drift to the left or to the right; every horizontal position in $\mathbb Z+\frac12$ is crossed by exactly two edges going from a vertex on the left hand-side site of this position to the right and two going from the right hand-side site to the left.
\end{itemize}
These rules leave six choices between neighbouring pairs of sites as shown below. We will refer to these as \emph{figures}.
\begin{figure}[h]
  \includegraphics[width=\linewidth]{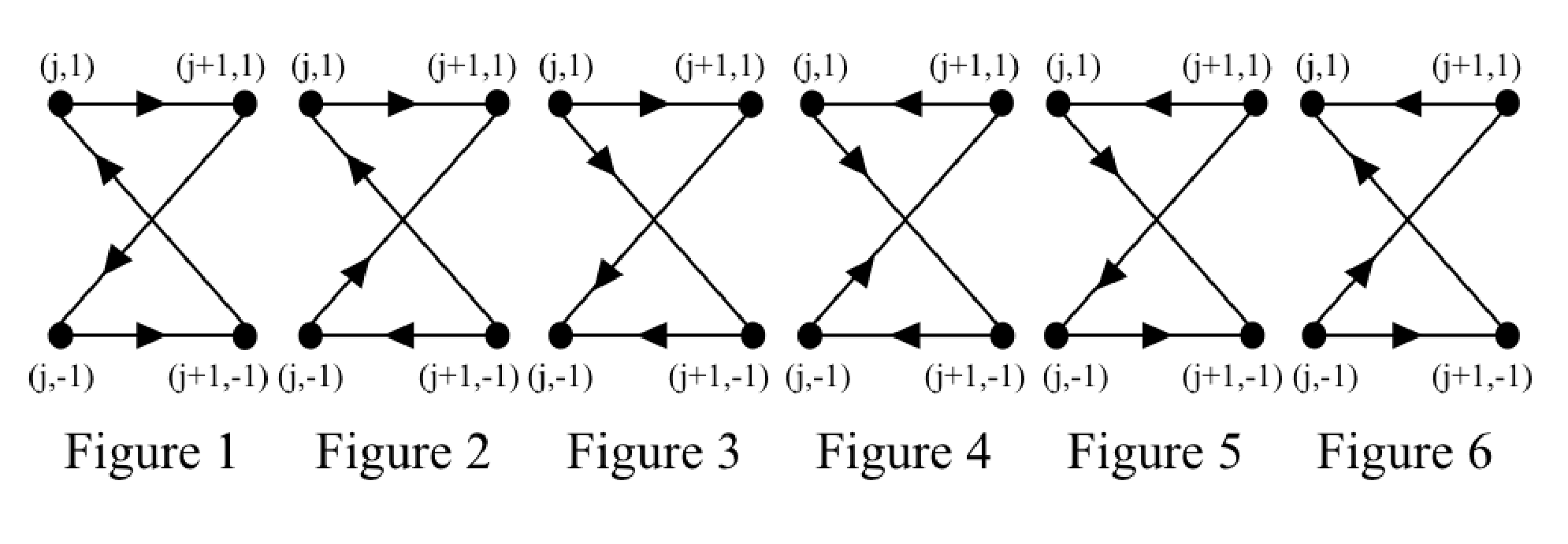}
\end{figure}
The divergence-free condition further restricts the orders in which the six figures can follow each other. Namely, each figure 2 must be immediately followed by figure 3, and each figure 5 must be immediately followed by figure 6. A figure 1, 3, 4 or 6 can be followed by any of the figures 1, 2, 4 or 5.

Without loss of generality, the graphs can be `twisted' vertically by replacing each figure 4 with figure 1 with no effect on the dynamics, since the only difference between figures 1 and 4 connecting two sites $j$ and $j+1$ is that the vertices $(j+1,1)$ and $(j+1,-1)$ are exchanged. However, in some cases where the number of figure 4's is odd it may be necessary to replace one of them with a figure 1 before `twisting' in order for the endpoints of the torus to match. Similarly we can replace each pair of figures 5 followed by 6 with a pair of figures 2 followed by 3. Therefore we only consider the case where the graph consists only of figures 1, 2 and 3.

Figures are thought of as being generated by a random process that obeys the rules above. We shall impose light ergodicity assumptions on this in Section \ref{sc:notndef}. It is this random orientation that our zero range process follows: nearest neighbour jumps on vertices are permitted across each edge in the direction of its orientation only. All through the paper we consider the quenched dynamics, i.e., the orientation of the edges is fixed and our results are valid for almost all realisations of this.

Next in Section \ref{sc:notndef} we introduce notation, make precise definitions and assumptions on the environments, and describe the local equilibrium measures. In Section \ref{sc:3} we state the main theorem, and in Section \ref{sc:4} we prove the hydrodynamic limit, following the steps of the relative entropy method in Kipnis and Landim \cite[pp. 115-130]{landim}. In Section \ref{sc:5} we prove a one block estimate required during the proof.

 \section{Notation and definitions}\label{sc:notndef}
 
 The following notation for the zero range process is taken from \cite[pp.28-30]{landim} and adapted where appropriate.
 
 The jump rate from a vertex with $k$ particles to an adjacent vertex to which jumps are permitted is given by $g(k)$, a nondecreasing function which satisfies $g(0) = 0.$ Let $\omega_x$ denote the number of particles at vertex $x$ and let the vector $\omega^{x,y} \in \Omega_N$ be defined by
 \[ (\omega^{x,y})_z = \begin{cases}
 	\omega_x - 1 & z = x \\
	\omega_y + 1 & z = y \\
	\omega_z & \text{otherwise.}
	\end{cases}.
\] 	
Letting $\text{type}(j)$ be the figure type between sites $j$ and $j+1$ and $\mathbb{N}_0 = \mathbb{N} \cup \{ 0 \}$, for any $N \in \mathbb{N}$ and function $f : \mathbb{N}_0^{\mathbb{T}_N \times \{ -1,1\}}$, the infinitesimal generator $L_N$ of the process is given by 
\[ \begin{aligned} L_N f(\omega) & = \sum_{\substack{j \in \mathbb{T}_N : \\ \text{type}(j) = 1}} \bigg{\{} g(\omega_{j,1})\Bigl( f(\omega^{(j,1),(j+1,1)})-f(\omega) \Bigr)  + g(\omega_{j,-1})\Bigl( f(\omega^{j,-1),(j+1,-1)})-f(\omega)\Bigr) \\ &\qquad
+ g(\omega_{j+1,1})\Bigl( f(\omega^{(j+1,1),(j,-1)})-f(\omega)\Bigr) + g(\omega_{j+1,-1})\Bigl( f(\omega^{(j+1,-1),(j,1)})-f(\omega) \Bigr)  \bigg{\}} \\ &\quad
+ \sum_{{\substack{j \in \mathbb{T}_N : \\ \text{type}(j) = 2}}} \bigg{\{} g(\omega_{j,1})\Bigl( f(\omega^{(j,1),(j+1,1)})-f(\omega) \Bigr)  + g(\omega_{j,-1})\Bigl( f(\omega^{(j,-1),(j+1,1)})-f(\omega)\Bigr) \\ & \qquad
+ g(\omega_{j+1,-1})\Bigl( f(\omega^{(j+1,-1),(j,1)})+f(\omega^{(j+1,-1),(j,-1)})-2f(\omega)\Bigr) \bigg{\}}  \\ & \quad
+ \sum_{{\substack{j \in \mathbb{T}_N : \\ \text{type}(j) = 3}}} \bigg{\{} g(\omega_{j,1}) \Bigl( f(\omega^{(j,1),(j+1,1)})+f(\omega^{(j,1),(j+1,-1)})-2f(\omega) \Bigr)  \\ & \qquad
+ g(\omega_{j+1,1})\Bigl( f(\omega^{(j+1,1),(j,-1)})-f(\omega)\Bigr) + g(\omega_{j+1,-1})\Bigl( f(\omega^{(j+1,-1),(j,-1)})-f(\omega) \Bigr) \bigg{\}}.
\end{aligned} \]  
  
 For any parameter $s > 0$ define the marginal probability measure $\hat{\nu}_s^1$ by
 \begin{equation*} \hat{\nu}_s^1(\omega_{(j,\pm 1)}=k)=\frac{s^k}{Z(s)g(k)!}  \end{equation*}
 where $Z(s)$ is a normalising constant that satisfies
 $$Z(s)=\sum_{k=0}^\infty \frac{s^k}{g(k)!},$$
and $g(k)!$ here means $\prod_{i=1}^k g(i)$ for $k \geq 1$ and $g(0)! = 1$. 
This measure is stationary for the zero range process, regardless of the value of $s$. Note that since $g$ is nondecreasing the power series $Z$ has a strictly positive radius of convergence $s^*$ at which the series diverges, and $\lim_{s \to s^*}Z(s) = \infty$. For a given density $\varrho$, we want $s$ to be such that $\mathbb{E}_{\hat{\nu}_s^1}[\omega_{(j, \pm 1)}]= \varrho$. Therefore $s$ must be a function of $\varrho$ and we say that $s = \Phi(\varrho)$ where 
 $\Phi$ is the inverse of
$$R(s)=\sum_{k=0}^\infty \frac{ks^k}{Z(s)g(k)!}.$$
 It is possible to choose such an $s$ because $R$ is a bijection from $[0,s^*)$ to $[0,\infty)$.

We set $\nu_\varrho^1 = \hat{\nu}_{\Phi(\varrho)}^1$ so that
 \begin{equation} \nu_\varrho^1(\omega_{(j,\pm 1)}=k) = \frac{{\Phi(\varrho)}^k}{Z(\Phi(\varrho))g(k)!} \label{eq:3} \end{equation}
 which ensures the expected number of particles at the vertices $(j,1)$ and $(j,-1)$ under the probability measure $\nu_\varrho^1$ is the density $\varrho$.
It follows that \begin{equation} \Phi(\varrho)=  \mathbb{E}_{\nu_{\varrho}^1}[g(k)]. \label{eq:phi} \end{equation} 
That is, the expected jump rate from a vertex with $k$ particles under the reference measure $\nu_\varrho^1$ is equal to the flux $\Phi$. 
  
The jump rate $g(k)$ is assumed to satisfy the condition 
$$ \mathbf{(SLG)} \ \limsup_{k \to \infty} \frac{g(k)}{k} = 0.$$
Kipnis and Landim \cite{landim}  proves the hydrodynamic limit of a mean-zero zero range process under both \textbf{(SLG)} and an alternative condition that the exponential moments of $g$ are finite.

For any $N \in \mathbb{N}$, define the measure
$$\nu^N=\otimes_{j \in \mathbb{T}_N} \Bigl( \nu_{1}^1 (\omega_{j,1}) \otimes \nu_{1}^1 (\omega_{j,-1}) \Bigr)$$
where the choice of 1 as a parameter is arbitrary and does not affect our calculations.
 It can easily be checked that $\nu^N$ is invariant. For a given density function $\rho : \mathbb{R}_{+} \times \mathbb{T} \to \mathbb{R}_{+}$ we also define the local equilibrium measure  
 $$\nu_{\rho(t,.)}^N (\omega)=\otimes_{j \in \mathbb{T}_N} \Bigl( \nu_{\rho(t,j/N)}^1(\omega_{j,1}) \otimes \nu_{\rho(t,j/N)}^1(\omega_{j,-1})\Bigr).$$ Both of these are product measures such that for $j \in \mathbb{T}_N$ both $\omega_{j,1}$ and $\omega_{j,-1}$ are distributed according to $\nu_{\rho(t,j/N)}^1$.

We introduce the local equilibrium measure $\tilde{\nu}_{\rho(t,.)}^N$ of second order approximation as in Funaki, Uchiyama and Yau, \cite{funaki} and Komoriya \cite{komoriya}.
This new measure depends on the specific figures in the graph, dividing it into tiles of various shapes within which the marginals at each site have the same distribution.
In particular, if a vertex has two inward arrows coming from the same direction it is included in the same tile that the inward arrows came from. 
However, if a vertex has one inward and one outward arrow in both directions then it belongs to a pair of vertices $(x,1)$ and $(x,-1)$ which share a tile. Every tile has exactly one such pair of vertices, which we call the `centre'. The tiles each  contain either two, three or four vertices.

This decomposition of the graph is illustrated below, with the figure numbers above the graph. 
\begin{figure}[h]
\includegraphics[width=\linewidth]{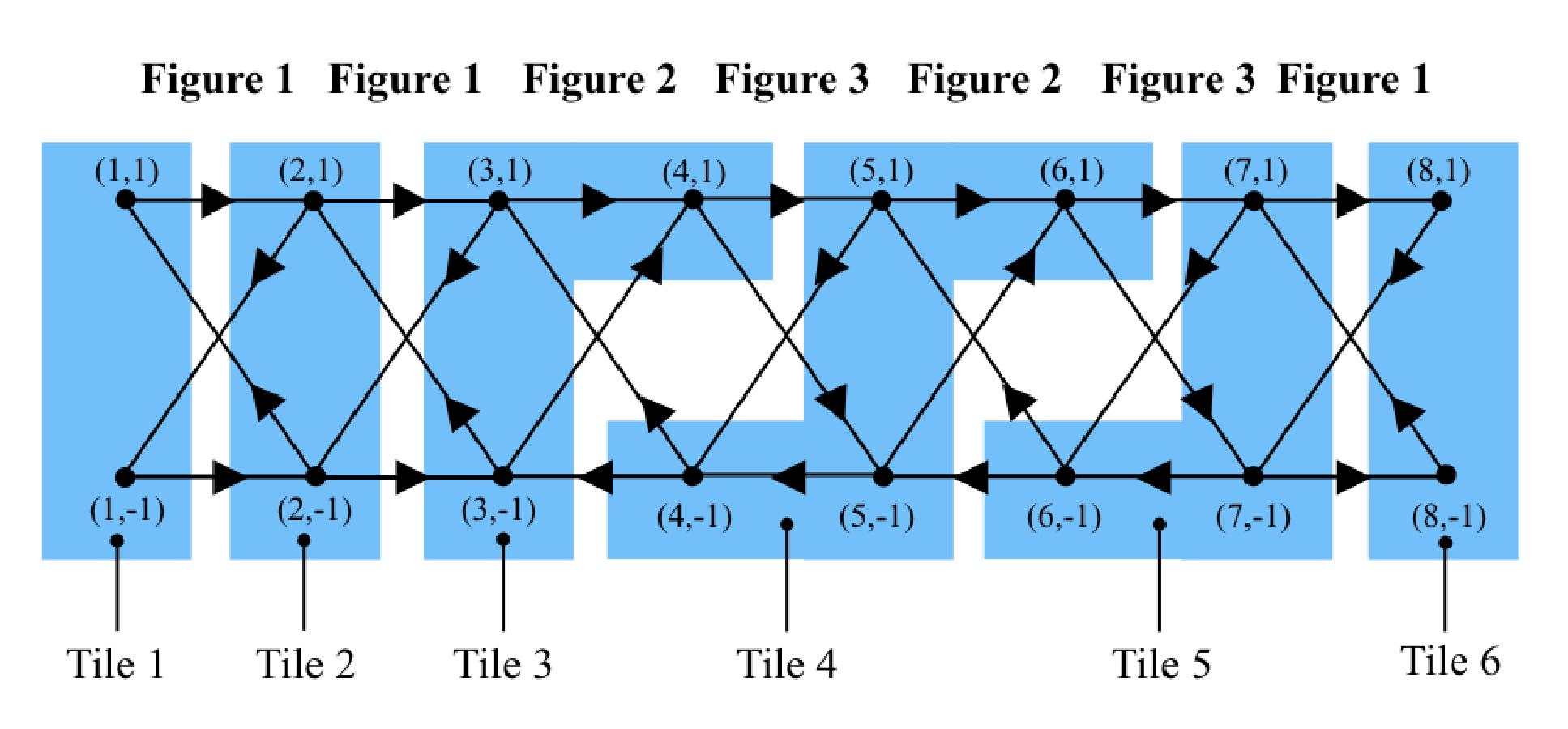}
\end{figure}
\newpage
For every graph there exists a unique decomposition into tiles according to the following algorithm, which uses the fact that each tile contains exactly one centre $(x,1)$ and $(x,-1)$, at most one of the two vertices $(x-1,1)$ and $(x-1,-1)$, and at most one of $(x+1,1)$ and $(x+1,-1)$.

\textbf{Algorithm.} Starting at vertices $(0,1)$ and $(0,-1)$, identify the two figures immediately to the left and right. These determine whether $(0,1)$ and $(0,-1)$ are in the same tile and, if so, the tile's shape. This is repeated for each pair of vertices $(x,1)$ and $(x,-1)$. The decisions made by the algorithm at each step are given below.

\begin{enumerate}
\item[$\ast$] If $\fbox{1} \begin{matrix} (x,1) \\ (x,-1) \end{matrix} \fbox{1}$  then $\begin{matrix} (x,1) \\ (x,-1) \end{matrix}$ is the centre of a tile containing $ (x,1)$, and $(x,-1)$ 
\item[$\ast$] If $\fbox{1} \begin{matrix} (x,1) \\ (x,-1) \end{matrix} \fbox{2}$  then $\begin{matrix} (x,1) \\ (x,-1) \end{matrix}$ is the centre of a tile containing $ (x,1)$, $(x,-1)$ and $(x+1,1)$ 
\item[$\ast$]  If $\fbox{2} \begin{matrix} (x,1) \\ (x,-1) \end{matrix} \fbox{3}$  then $\begin{matrix} (x,1) \\ (x,-1) \end{matrix}$ is not the centre of a tile
\item[$\ast$] If $\fbox{3} \begin{matrix} (x,1) \\ (x,-1) \end{matrix} \fbox{1}$  then $\begin{matrix} (x,1) \\ (x,-1) \end{matrix}$ is the centre of a tile containing $(x-1,-1)$, $ (x,1)$ and $(x,-1)$
\item[$\ast$] If $\fbox{3} \begin{matrix} (x,1) \\ (x,-1) \end{matrix} \fbox{2}$  then $\begin{matrix} (x,1) \\ (x,-1) \end{matrix}$ is the centre of a tile containing $(x-1,-1)$, $ (x,1)$, $(x,-1)$ and $(x+1,1)$
\end{enumerate}

To generate these graphs on the torus, we fix a common random sequence $\sigma \in \{ 1,2,3 \}^{\mathbb{N}}$ according to some probability measure $\mathcal{P}$, subject to the following conditions.
\begin{itemize}
 \item[\textbf{(G1)}] The figures generated are ergodic with respect to the figure shift operator.
 \item[\textbf{(G2)}] For any $l > 0$ and any sequence of $l$ tiles which is possible according to the above algorithm, the number of times this sequence occurs in the first \(N\) terms tends to infinity as $N \to \infty$.
\end{itemize}

Such a sequence could, for example, arise via a Markov chain progressing through $\mathbb N$ and at each point where a choice can be made (that is, after a `1' 
 or at the end of a `2-3'
 pair), starting a pair with probability $p$ and inserting a `1'
 otherwise.

The graph on the torus \(\mathbb{T}_N\) is then drawn as the first $N$ terms of $\sigma$ and the final terms are altered if necessary for the two endpoints of the torus to be compatible. If the sequence begins with figure `3' and the $(N-1)$-th figure is `2', we change figure $N-1$ to `1' and figure $N$ to `2', since there are no figures that can be placed between `2' and `3'. Otherwise, if the endpoints are not compatible we select the $N$-th figure uniformly at random from the figures that make the endpoints compatible. Since we only alter at most two figures, we will be able to make advantage of the conditions $\textbf{(G1)}$ and $\textbf{(G2)}$ that were valid for the infinite generating sequence. This construction also provides a common realisation of graphs on \(\mathbb T_N\) across different \(N\) values, hence a.s.\ or in probability statements as \(N\to\infty\) will make sense.

We enumerate tiles according to the site of their centre, so that the first tile has centre $1 \in \mathbb{T}_N$ unless the site 1 is not the centre of a tile, in which case the first tile has centre $2 \in \mathbb{T}_N$. 
 Let $T_N$ denote the number of tiles in the $N$-th graph of the sequence. 
The ergodicity condition implies the existence of the constant
 \begin{equation} \kappa := \lim_{N \to \infty} \frac{N}{T_N}. \label{eq:kappa} \end{equation}
  We note that $\kappa$ is a deterministic constant that exists almost surely, and that altering at most two figures between the generating sequence and the one on \(\mathbb T_N\) for each graph does not affect the limit.
 
 The measure $\tilde{\nu}_{\rho(t,.)}^N$ is defined as the product of marginals given by \eqref{eq:3} with parameter $\rho(t,j/T_N)$ for vertices in tile $j$ of the graph.

Let $x_j \in \mathbb{T}_N$ denote the site containing the unique pair of vertices $(x_j,1)$ and $(x_j,-1)$ which are both contained in tile $j$ of the graph.

Given a measure $\mu^N$ corresponding to an initial configuration, let $\mu_t^N=S_{tN^2}\mu^N$ where $S_t$ is the semigroup associated with $L_N$. The Radon-Nikodym derivatives are given by
$$f_t^N = \frac{\text{d}\mu_t^N}{\text{d}\nu^N}, \quad \psi_t^N = \frac{\text{d}\nu_{\rho(t,.)}^N}{\text{d}\nu^N}, \quad \tilde{\psi}_t^N = \frac{\text{d}\tilde{\nu}_{\rho(t,.)}^N}{\text{d}\nu^N}.$$
The relative entropy with respect to the usual local equilibrium measure $\nu_{\rho(s,.)}^N$ is defined by
$$H_N(t) = H(\mu_t^N | \nu_{\rho(t,.)}^N) = \int_{\Omega_N} f_t^N(\omega)\log \frac{f_t^N(\omega)}{\psi_t^N(\omega)}\text{d}\nu^N(\omega),$$
and the relative entropy with respect to $\tilde{\nu}_{\rho(s,.)}^N$ is defined by
$$\tilde{H}_N(t) = \tilde{H}(\mu_t^N | \tilde{\nu}_{\rho(t,.)}^N) = \int_{\Omega_N} f_t^N(\omega)\log \frac{f_t^N(\omega)}{\tilde{\psi}_t^N(\omega)}\text{d}\nu^N(\omega).$$

For any parameter $\rho(s,u)$ let $(X_j)_{j=1}^\infty$ be a sequence of independent, identically distributed random variables with the same distribution as a marginal of $\nu_{\rho(s,.)}^N$ at a site with density $\rho(s,u)$.  That is, \begin{equation} \mathbb{P}_{\rho(s,u)}(X_j=k) = \nu_{\rho(s,u)}^1(\omega=k). \label{eq:4}
\end{equation}
 Let $\mathbb{E}_{\rho(s,u)}$ be the expectation under the probability distribution $\mathbb{P}_{\rho(s,u)}$.

\section{Results}\label{sc:3}

With particle density $\rho(t,x)$ at the macroscopic time $t$ under diffusive scaling and macroscopic horizontal position $x \in \mathbb{T} = \mathbb{R}/ \mathbb{Z}$, the hydrodynamic limit of this system will be \begin{equation}
\partial_t \rho(t,x)=\kappa \partial_x^2 \Phi(\rho(t,x)), \ (t,x) \in (0,\infty)\times \mathbb{T} \label{eq:1}
\end{equation}
\begin{equation}
\rho(0,x)=\rho_0(x), \ x \in \mathbb{T} \label{eq:2}
\end{equation}
where $\rho_0$ is a twice continuously differentiable initial condition, $\Phi$ is given by \eqref{eq:phi} and $\kappa$ by \eqref{eq:kappa}. The solution of \eqref{eq:1} and \eqref{eq:2} is smooth for all $t \geq 0$. The proof of the hydrodynamic limit requires $\rho(t,x)$ to be bounded both above and away from zero, so we impose the condition $K_1 \leq \rho_0(x) \leq K_2$ for all $x \in \mathbb{T}$ for some $K_1, K_2 > 0$. This ensures $K_1 \leq \rho(t,x) \leq K_2$ for all $t> 0$ due to the maximum principle since $\Phi$ is monotone increasing. In the case where $\kappa = 1$ this is the same PDE as obtained in \cite{landim} for the symmetric mean-zero zero range process on $\mathbb{T}_N$.

The following theorem, which states that the local particle distribution converges to the solution of \eqref{eq:1} and \eqref{eq:2} in probability, will be proved in Sections \ref{sc:4.2} to \ref{sc:4.5} using the relative entropy method. Here, $H_N(t)$ denotes the relative entropy of the local equilibrium measure $\nu_{\rho(t,.)}^N$ with respect to $\mu_t^N$.

\begin{theorem}\label{tm:1} Assume \textbf{(G1)} and \textbf{(G2)} on the environment and that the jump rate function $g$ is nondecreasing and satisfies \textbf{(SLG)} and $g(0) = 0$. Let $\rho$ be a solution of \eqref{eq:1} and \eqref{eq:2}.
Then 
for $\mathcal{P}$-almost any random environment $\sigma \in \{ 1,2,3 \}^{\mathbb{N}}$,
 if $H_N(0) = o(N)$ then the density of particles in $\mathbb{T}_N \times \{-1,1\}$ converges in probability to $\rho(t,x)\text{d}x$ for all $t$. That is, letting $\phi \in C^{\infty}(\mathbb{T})$ and letting $\psi:\mathbb{Z}^{\{1,-1\}\times \mathbb{Z}} \to \mathbb{R}$ be a bounded function, then
$$\lim_{N \to \infty}\mathbb{E}_{\mu_t^N}\Big{[}\Big{|} \frac{1}{N}\sum_{j \in \mathbb{T}_N}\phi(j/N)\psi(\tau_j \omega)-\int_{\mathbb{T}}\phi(x)\mathbb{E}_{\nu_{\rho(t,x)}}(\psi(\omega))\text{d}x \Big{|}\Big{]} = 0.$$
\end{theorem}

\section{Relative entropy method}\label{sc:4}

\subsection{Outline of proof}

  The first step in proving Theorem \ref{tm:1} is to show that if $H_N(0) = o(N)$ then $\tilde{H}_N(t) = o(N)$ for any $t > 0$. This is achieved by proving in Lemma \ref{lm:5} that $H_N(0) = o(N)$ implies $\tilde{H}_N(0)=o(N)$, and then bounding $\partial_t \tilde{H}_N(t)$ to show that 
 $$\tilde{H}_N(t) = \tilde{H}_N(0) + \int_0^t \partial_s \tilde{H}_N(s) \ \text{d}s \leq K \int_{0}^t \tilde{H}_N(s)\text{d}s + o(N)$$ for some constant $K > 0$ to be determined later. We then use Gronwall's inequality to conclude that $\tilde{H}_N(t) = o(N)$.
 It is necessary to work with the measure $\tilde{\nu}_{\rho(t,.)}^N$ and show that $\tilde{H}_N(t) = o(N)$ instead of $H_N(t) = o(N)$ since $\tilde{\nu}_{\rho(t,.)}^N$ divides the graph into tiles in such a way that non-gradient terms do not arise when computing the adjoint of the infinitesimal generator. This is because the only contributions to the expression are from vertices containing inward and outward arrows from both directions that allow a second order Taylor expansion to be performed.

  The result that $H_N(0) = o(N)$ implies $\tilde{H}_N(t) = o(N)$ obtained in Proposition \ref{pr:6} allows us to prove Theorem \ref{tm:1}, analogously to the method used in Kipnis and Landim \cite{landim}.

\subsection{Bounding $\partial_t \tilde{H}_N(t)$ - Initial Steps}\label{sc:4.2}

As in \cite[p.9]{funaki}, the time derivative of the relative entropy satisfies the inequality $$\partial_t \tilde{H}_N(t) \leq  \int_{\Omega_N} \frac{1}{\tilde{\psi}_t^N}(N^2 L_N^* \tilde{\psi}_t^N - \partial_t \tilde{\psi}_t^N )f_t^N \text{d}\nu^N,$$
where $L_N^*$ is the adjoint of $L_N$. Note that $L_N^*$ is the generator of the time-reversed process for which the arrows are reversed and the dynamics are otherwise the same.

Where the site $s \in \mathbb{T}_N$ is the centre of a tile, let $k_s$ denote the corresponding number of that tile. Due to the product structure of the measure $\tilde{\nu}_{\rho(t,.)}^N$, the terms of $\frac{L_N^* \tilde{\psi}_t^N}{\tilde{\psi}_t^N}$ have the form  
\[ \begin{aligned} & 
\frac{g(\omega_{i,1})}{\tilde{\psi}_t^N(\omega)} \Bigl( \tilde{\psi}_t^N(\omega^{(i,1),(i+1,-1)})-\tilde{\psi}_t^N(\omega)\Bigr) \\
& = g(\omega_{i,1})\Bigg{\{} \frac{\Phi \Bigl( \rho \Bigl( t,\frac{k_i}{T_N} \Bigr) \Bigr)^{\omega_{i,1}-1} \Phi \Bigl( \rho \Bigl( t,\frac{k_{i+1}}{T_N}\Bigr) \Bigr)^{\omega_{i+1,-1}+1}}{\Phi \Bigl( \rho \Bigl( t,\frac{k_i}{T_N} \Bigr) \Bigr)^{\omega_{i,1}} \Phi \Bigl( \rho \Bigl( t,\frac{k_{i+1}}{T_N} \Bigr) \Bigr)^{\omega_{i+1,-1}}} - 1 \Bigg{\}} \\
& = \frac{g(\omega_{i,1})}{\Phi \Bigl( \rho \Bigl( t,\frac{k_i}{T_N} \Bigr) \Bigr) }\Bigg{\{} \Phi \Bigl( \rho \Bigl( t,\frac{k_{i+1}}{T_N}\Bigr)\Bigr) - \Phi \Bigl( \rho \Bigl( t,\frac{k_i}{T_N} \Bigr) \Bigr) 
\Bigg{\}}
\end{aligned} \]
where, in this example, there is a figure 1 connecting sites $i$ and $i+1$.
  We now compute $\frac{N^2 L_N^* \tilde{\psi}_t^N}{\tilde{\psi}_t^N}$ and $\frac{\partial_t \tilde{\psi}_t^N}{\tilde{\psi}_t^N}$.
\begin{equation} \begin{aligned} & \frac{N^2L_N^*\tilde{\psi}_t^N}{\tilde{\psi}_t^N}  
\\ & = N^2 \sum_{\substack{i \in \mathbb{T}_N : \\ \text{type}(i-1) = 1 \\ \text{type}(i) = 1}} \frac{g(\omega_{i,1})+g(\omega_{i,-1})}{\Phi \Bigl( \rho \Bigl( t,\frac{k_{i}}{T_N} \Bigr) \Bigr)}  \bigg{\{}   \Phi \Bigl( \rho \Bigl( t,\frac{k_{i-1}}{T_N} \Bigr) \Bigr) - 2\Phi \Bigl( \rho \Bigl( t,\frac{k_{i}}{T_N} \Bigr) \Bigr)    + \Phi \Bigl( \rho \Bigl( t,\frac{k_{i+1}}{T_N} \Bigr) \Bigr) \bigg{\}}
\\ & \quad + N^2 \sum_{\substack{i \in \mathbb{T}_N : \\ \text{type}(i-1) = 1 \\ \text{type}(i) = 2}} \frac{g(\omega_{i,1}) + g(\omega_{i,-1})}{\Phi \Bigl( \rho \Bigl( t,\frac{k_{i}}{T_N} \Bigr) \Bigr)} \bigg{\{}   \Phi \Bigl( \rho \Bigl( t,\frac{k_{i-1}}{T_N} \Bigr) \Bigr) - 2\Phi \Bigl( \rho \Bigl( t,\frac{k_{i}}{T_N} \Bigr) \Bigr)  + \Phi \Bigl( \rho \Bigl( t,\frac{k_{i+2}}{T_N} \Bigr) \Bigr) \bigg{\}}
\\ & \quad + N^2 \sum_{\substack{i \in \mathbb{T}_N : \\ \text{type}(i-1) = 2 \\ \text{type}(i) = 3}} \Bigg{[} \frac{g(\omega_{i,1})}{\Phi \Bigl( \rho \Bigl( t,\frac{k_{i-1}}{T_N} \Bigr) \Bigr)} \bigg{\{}   2\Phi \Bigl( \rho \Bigl( t,\frac{k_{i-1}}{T_N} \Bigr) \Bigr) - 2\Phi \Bigl( \rho \Bigl( t,\frac{k_{i-1}}{T_N} \Bigr) \Bigr) \bigg{\}} \\ &  \qquad + \frac{g(\omega_{i,-1})}{\Phi \Bigl( \rho \Bigl( t,\frac{k_{i-1}}{T_N} \Bigr) \Bigr)} \bigg{\{}   2\Phi \Bigl( \rho \Bigl( t,\frac{k_{i+1}}{T_N} \Bigr) \Bigr) - 2\Phi \Bigl( \rho \Bigl( t,\frac{k_{i+1}}{T_N} \Bigr) \Bigr) \bigg{\}} \Bigg{]}
\\ & \quad + N^2 \sum_{\substack{i \in \mathbb{T}_N : \\ \text{type}(i-1) = 3 \\ \text{type}(i) = 2}} \frac{g(\omega_{i,1}) + g(\omega_{i,-1})}{\Phi \Bigl( \rho \Bigl( t,\frac{k_{i}}{T_N} \Bigr) \Bigr)} \bigg{\{}   \Phi \Bigl( \rho \Bigl( t,\frac{k_{i-2}}{T_N} \Bigr) \Bigr) - 2\Phi \Bigl( \rho \Bigl( t,\frac{k_{i}}{T_N} \Bigr) \Bigr)   + \Phi \Bigl( \rho \Bigl( t,\frac{k_{i+2}}{T_N} \Bigr) \Bigr) \bigg{\}}
\\ & \quad + N^2 \sum_{\substack{i \in \mathbb{T}_N : \\ \text{type}(i-1) = 3 \\ \text{type}(i) = 1}} \frac{g(\omega_{i,1}) + g(\omega_{i,-1})}{\Phi \Bigl( \rho \Bigl( t,\frac{k_{i}}{T_N} \Bigr) \Bigr)} \bigg{\{}   \Phi \Bigl( \rho \Bigl( t,\frac{k_{i-2}}{T_N} \Bigr) \Bigr) - 2\Phi \Bigl( \rho \Bigl( t,\frac{k_{i}}{T_N} \Bigr) \Bigr)   + \Phi \Bigl( \rho \Bigl( t,\frac{k_{i+1}}{T_N} \Bigr) \Bigr) \bigg{\}}
\\ &  = N^2 \sum_{j=1}^{T_N} \frac{g(\omega_{x_j,1})+g(\omega_{x_j,-1})}{\Phi \Bigl( \rho \Bigl( t,\frac{j}{T_N} \Bigr) \Bigr) }\bigg{\{} \Phi \Bigl( \rho \Bigl( t,\frac{j-1}{T_N} \Bigr) \Bigr) -2\Phi \Bigl( \rho \Bigl( t,\frac{j}{T_N} \Bigr) \Bigr)+\Phi \Bigl( \rho \Bigl( t,\frac{j+1}{T_N} \Bigr) \Bigr) \bigg{\}}  \\ & = 
\frac{N^2}{T_N^2} \sum_{j=1}^{T_N} (g(\omega_{x_j,1})+g(\omega_{x_j,-1}))\frac{\partial_x^2 \Phi(\rho(t,j/T_N))}{\Phi(\rho(t,j/T_N))} + o(N) \end{aligned} \label{eq:5} \end{equation}

using the Taylor expansion of $\Phi(\rho(t,\frac{j}{T_N}))$ in the last line. Observe that the third sum, corresponding to sites $i \in \mathbb{T}_N$ between a figure 2 and a figure 3, vanishes so that there are no non-gradient terms.
The factor of $\frac{N^2}{T_N^2}$ in the last line is due to the step size.
The error is $o(N)$ because we assume $\rho_0$ is smooth enough that its second derivative is H\"older continuous and thus $\partial_x^2 \rho(t,x)$ is H\"older continuous in $x$ and the remainder term in the Taylor expansion is bounded by a H\"older norm.

We can subtract the telescopic expression 
\begin{equation} \begin{aligned} 0 &=  2N^2 \sum_{j=1}^{T_N} \Phi(\rho(t,j/T_N))\Bigl( \frac{\Phi(\rho(t,(j+1)/T_N))}{\Phi(\rho(t,j/T_N))}+\frac{\Phi(\rho(t,(j-1)/T_N))}{\Phi(\rho(t,j/T_N))}-2 \Bigr) \\ & = \frac{2N^2}{T_N^2} \sum_{j=1}^{T_N} \Phi(\rho(t,j/T_N)) \frac{\partial_x^2 \Phi(\rho(t,j/T_N))}{\Phi(\rho(t,j/T_N))} + o(N) \end{aligned} \label{eq:6} \end{equation}
which is equal to zero. The first and second lines are equal due to the same Taylor expansion as in
\eqref{eq:5}.

Let $N_j$ be the number of vertices in tile $j$ of the graph and let $\hat{\omega}_j$ be the total number of particles on  tile $j$. Then the time derivative is given by
\begin{equation} \begin{aligned} \frac{\partial_t \tilde{\psi}_t^N}{\tilde{\psi}_t^N} & =\partial_t \log \tilde{\psi}_t^N =\sum_{j=1}^{T_N} \Bigl(  \hat{\omega}_j \frac{\Phi'(\rho(t,j/T_N))\partial_t \rho(t,j/T_N)}{\Phi(\rho(t,j/T_N))} \\ & \quad -N_j \frac{Z'(\Phi(\rho(t,j/T_N)))\Phi'(\rho(t,j/T_N))\partial_t \rho(t,j/T_N)}{Z(\Phi(\rho(t,j/T_N)))}\Bigr) \\ &
= \sum_{j=1}^{T_N}  \Bigl( \hat{\omega}_j-N_j\rho(t,j/T_N) \Bigr) \frac{\kappa \Phi'(\rho(t,j/T_N))\partial_x^2 \Phi(\rho(t,j/T_N))}{\Phi(\rho(t,j/T_N))} \label{eq:7} \end{aligned}
 \end{equation}
 where in the last line we have substituted in the PDE $\partial_t \rho(t,x)=\kappa \partial_x^2 \Phi(\rho(t,x))$ and used the fact that $$\frac{Z'(\Phi(\rho(t,x)))}{Z(\Phi(\rho(t,x)))} = \frac{R(\Phi(\rho(t,x)))}{\Phi(\rho(t,x))} = \frac{\rho(t,x)}{\Phi(\rho(t,x))}.$$

Note that the integral of $\cfrac{1}{\tilde{\psi}_t^N}(N^2 L_N^* \tilde{\psi}_t^N - \partial_t \tilde{\psi}_t^N)$ is the same as the integral of the sum of \eqref{eq:5}, \eqref{eq:6} and \eqref{eq:7} plus an $o(N)$ error. That is,
\begin{equation}
 \begin{aligned}
  &\int_{\Omega_N}\frac{1}{\tilde{\psi}_t^N}(N^2 L_N^* \tilde{\psi}_t^N - \partial_t \tilde{\psi}_t^N )f_t^N \text{d}\nu^N\\
  &=\int_{\Omega_N} \sum_{j=1}^{T_N} \Bigg{\{} \frac{N^2}{T_N^2} (g(\omega_{x_j,1})+g(\omega_{x_j,-1}))\frac{\partial_x^2 \Phi(\rho(t,j/T_N))}{\Phi(\rho(t,j/T_N))}
 -\frac{2N^2}{T_N^2} \Phi(\rho(t,j/T_N)) \frac{\partial_x^2 \Phi(\rho(t,j/T_N))}{\Phi(\rho(t,j/T_N))}\\
  &\quad - \Bigl( \hat{\omega}_j-N_j\rho(t,j/T_N) \Bigr) \frac{\kappa \Phi'(\rho(t,j/T_N))\partial_x^2 \Phi(\rho(t,j/T_N))}{\Phi(\rho(t,j/T_N))} \Bigg{\}} f_t^N \text{d}\nu^N + o(N).
 \end{aligned}\label{eq:beta}
\end{equation}

\subsection{Replacements and application of the one block estimate}\label{sc:4.3}

For the $2l+1$ tiles surrounding the $j$-th tile, define the average particle density per vertex $$\omega_j^l = \frac{1}{\sum_{|k-j| \leq l}N_k}\sum_{|k-j| \leq l} \hat{\omega}_k$$ and the average particle density per tile $$\bar{\omega}_j^l = \frac{1}{2l+1}\sum_{| k - j | \leq l} \hat{\omega}_k.$$ Note that for large values of $l$, the typical values of $\frac{1}{\kappa} \bar{\omega}_j^l$ are close to $2 \omega_j^l$.

In this section we will perform replacements which enable us to write that 
\begin{equation} \begin{aligned}
 &\int_{\Omega_N} \cfrac{1}{\tilde{\psi}_t^N}(N^2 L_N^* \tilde{\psi}_t^N - \partial_t \tilde{\psi}_t^N)  f_t^N \text{d}\nu^N \\ & =
\int_{\Omega_N} \sum_{j=1}^{T_N} 2 \kappa F(t,j/T_N) \Big{\{} \Phi(\omega_j^l) - \Phi(\rho(t,j/T_N)) - \Phi'(\rho(t,j/T_N))( \omega_j^l - \rho(t,j/T_N)) \Big{\}}f_t^N \text{d}\nu^N \\ & \quad + o(N) + O(N)C(l) \end{aligned} \label{eq:8} \end{equation}
where $$F(t,j/T_N) = \frac{\kappa \partial_x^2 \Phi(\rho(t,j/T_N))}{\Phi(\rho(t,j/T_N))}$$
and $C(l)$ is such that $\lim_{l \to \infty}C(l) = 0$. 
In particular, recalling that $\frac{N}{T_N} \to \kappa$ as $N \to \infty$ in \eqref{eq:5} and \eqref{eq:6}, we will replace each $\hat{\omega}_j$ in 
\eqref{eq:7} with $\bar{\omega}_j^l$ and then $2\kappa \omega_j^l$, and replace each $N_j$ factor  with $2\kappa$. 

 Replacing $\hat{\omega}_j$ with $\bar{\omega}_j^l$ in the last line of \eqref{eq:beta} creates the error term
\begin{equation} \begin{aligned} 
\sum_{j=1}^{T_N} &  \Bigl( \hat{\omega}_{j} - \frac{1}{2l+1}\sum_{|k-j|\leq l}\hat{\omega}_k \Bigr)F(t,j/T_N)\Phi'(\rho(t,j/T_N))
\\ & = \frac{1}{2l+1}\sum_{j=1}^{T_N} \hat{\omega}_j \Bigl( (2l+1)F(t,j/T_N)\Phi'(\rho(t,j/T_N))-\sum_{|k-j| \leq l} F(t,k/T_N)\Phi'(\rho(t,k/T_N)) \Bigr). \end{aligned} \label{eq:9} \end{equation}
Performing first order Taylor expansions to extract factors of $\frac{l}{T_N}$ and using the fact that $T_N \geq N/2$,  each $$ \hat{\omega}_j \Bigl( (2l+1)F(t,j/T_N)\Phi'(\rho(t,j/T_N))-\sum_{|k-j| \leq l} F(t,k/T_N)\Phi'(\rho(t,k/T_N)) \Bigr) $$ expression is bounded by $$ \frac{2l (2l+1) }{N} \Big{\|} \partial_x(F(t,x)\Phi'(\rho(t,x))) \Big{\|}_\infty  \hat{\omega}_j ,$$ so that the integral of the error term \eqref{eq:9} with respect to $\mu_t^N $ is bounded by
$$ \frac{2l(2l+1)}{N} \Big{\|} \partial_x(F(t,x)\Phi'(\rho(t,x))) \Big{\|}_\infty \mathbb{E}_{\mu_t^N} \Big{[} \sum_{j=1}^{T_N} \hat{\omega}_j \Big{]}.$$
To show that this expression is $o(N)$, first observe that the range of $\rho$ is bounded due to the initial condition \eqref{eq:2} and on this range $\Phi'$ and thus  $ \| \partial_x (F(t,x)\Phi'(\rho(t,x))) \|_\infty$ is bounded.
 Furthermore $T_N \geq N/2$, and the following lemma can be applied.

\begin{lemma}\label{lm:2}
 $\frac{1}{N} \mathbb{E}_{\mu_t^N} \Big{[} \sum_{j=1}^{T_N} \hat{\omega}_j \Big{]}$ is bounded.
\end{lemma}

\begin{proof} 
This proof is adapted from Kipnis and Landim \cite[p.84]{landim}. Since the system is conservative, the total number of particles at time $t$ is the same as time $0$. Applying the entropy inequality at time 0, the expression is bounded by
$$\frac{1}{\gamma N}\tilde{H}_N(0) + \frac{1}{\gamma N} \sum_{j=1}^{T_N} \log \mathbb{E}_{\tilde{\nu}_{\rho(0,.)}^N} [e^{\gamma \hat{\omega}_j}].$$
It will be shown later that $\tilde{H}_N(0) = o(N)$, which leaves the second term. By independence, the expression inside each logarithm is a product of $\mathbb{E}_{\tilde{\nu}_{\rho(0,.)}^N}[e^{\gamma \omega_x}]$ for at most four vertices $x$, and the Laplace transform of $\omega_x$ is finite if $\gamma$ is chosen to be sufficiently small.
\end{proof}

Next we will replace $ \bar{\omega}_j^l$ by $2 \kappa \omega_j^l$. The error term resulting from this replacement is

\begin{equation*} \begin{aligned} & \mathbb{E}_{\mu_t^N} \bigg{[} \sum_{j=1}^{T_N}  \kappa F(t,j/T_N)\Phi'(\rho(t,j/T_N)) \bigg{(} \frac{1}{\kappa} - \frac{4l+2}{\sum_{|k - j| \leq l} N_k} \bigg{)} \biggl( \frac{1}{2l+1} \sum_{|k - j | \leq l} \hat{\omega}_k \biggr) \bigg{]} \\ &  =  \mathbb{E}_{\mu_t^N}  \bigg{[} \sum_{j=1}^{T_N}  \frac{\kappa \hat{\omega}_j}{2l+1}\sum_{|k-j| \leq l} F(t,k/T_N)\Phi'(\rho(t,k/T_N)) \Big{(} \frac{1}{\kappa }-\frac{4l+2}{\sum_{|m-k| \leq l} N_m}  \Big{)}  \bigg{]} \\ &
\leq  \mathbb{E}_{\mu_t^N}   \bigg{[}  \Big{\|} \kappa F(t,x)\Phi'(\rho(t,x)) \Big{\|}_\infty \sum_{j=1}^{T_N} \frac{\hat{\omega}_j}{2l+1} \sum_{|k - j| \leq l} \Big{|} \frac{1}{\kappa }-\frac{4l+2}{\sum_{|m-k| \leq l} N_m}  \Big{|}   \bigg{]}. \end{aligned}
 \end{equation*}
Applying the entropy inequality,
\begin{equation*}
 \begin{aligned}
  & \mathbb{E}_{\mu_t^N}   \bigg{[} \Big{\|} \kappa F(t,x)\Phi'(\rho(t,x)) \Big{\|}_\infty \sum_{j=1}^{T_N} \frac{\hat{\omega}_j}{2l+1} \sum_{|k - j| \leq l} \Big{|} \frac{1}{\kappa }-\frac{4l+2}{\sum_{|m-k| \leq l} N_m}  \Big{|}   \bigg{]} \\
  & \quad
\leq \frac{1}{\gamma}\tilde{H}_N(t)
 +\frac{1}{\gamma} \sum_{j=1}^{T_N} \log \mathbb{E}_{\tilde{\nu}_{\rho(t,.)}^N} \bigg{[} \exp \Bigl( \frac{\| \kappa F(t,x)\Phi'(\rho(t,x)) \|_\infty \gamma \hat{\omega}_j}{2l+1} \sum_{|k - j| \leq l} \Big{|} \frac{1}{\kappa }-\frac{4l+2}{\sum_{|m-k| \leq l} N_m}  \Big{|} \Bigr) \bigg{]}. \end{aligned}
\end{equation*}
For $\gamma$ sufficiently small each Laplace transform
$$\mathbb{E}_{\tilde{\nu}_{\rho(t,.)}^N} [e^{\| \kappa F(t,x)\Phi'(\rho(t,x)) \|_\infty \gamma \hat{\omega}_j}]$$
 is bounded, and the $\frac{1}{2l+1}\sum_{|k-j| \leq l} \Big{|} \frac{1}{\kappa }-\frac{4l+2}{\sum_{|m-k| \leq l} N_m}  \Big{|}$ factors are as well. Applying Jensen's inequality in the form $$\mathbb{E}[f(\omega)^p] \leq \mathbb{E}[f(\omega)]^p$$ with exponent $$p = \frac{\sum_{|k-j| \leq l} \Big{|} \frac{1}{\kappa }-\frac{4l+2}{\sum_{|m-k| \leq l} N_m}  \Big{|}}{2l+1} < 1,$$ which is less than 1 because the terms $\frac{1}{\kappa}$ and $\frac{4l+2}{\sum_{|m-k| \leq l}N_m}$ inside the modulus sign both lie between 0 and 1, the second term of the right hand side of the inequality is bounded by
$$\frac{1}{\gamma}\sum_{j=1}^{T_N} \Bigl( \sum_{|k-j| \leq l} \frac{1}{2l+1} \bigg{|} \frac{1}{\kappa }-\frac{4l+2}{\sum_{|m-k| \leq l} N_m}  \bigg{|} \Bigr) \log \mathbb{E}_{\tilde{\nu}_{\rho(t,.)}^N} [e^{\| \kappa F(t,x)\Phi'(\rho(t,x)) \|_\infty \gamma \hat{\omega}_j}].$$
This can be rewritten as
$$\frac{1}{\gamma}\sum_{j=1}^{T_N} \bigg{|} \frac{1}{\kappa }-\frac{4l+2}{\sum_{|m-j| \leq l} N_m}  \bigg{|} \log \mathbb{E}_{\tilde{\nu}_{\rho(t,.)}^N} [e^{\| \kappa F(t,x)\Phi'(\rho(t,x)) \|_\infty \gamma \hat{\omega}_j}].$$
Each expectation is bounded by a constant. Since we are aiming to show this error term is $o(N)$, we divide by $N$.  
  
Let $\langle . \rangle$ denote the ergodic limit of an expression as $N \to \infty$.
 Due to the ergodicity condition with respect to the figure shift operator \textbf{(G1)}, the sequence of $$\bigg{|} \frac{1}{\kappa }-\frac{4l+2}{\sum_{|m-j| \leq l} N_m}  \bigg{|}$$ values converges to its ergodic limit as $N \to \infty$. Hence
 $$\frac{1}{N}\sum_{j=1}^{T_N} \bigg{|} \frac{1}{\kappa }-\frac{4l+2}{\sum_{|m-j| \leq l} N_m}  \bigg{|} \to \lim_{N \to \infty}\frac{T_N}{N} \bigg{\langle} \bigg{|} \frac{1}{\kappa} - \frac{4l+2}{\sum_{j=1}^{2l+1}N_j}\bigg{|} \bigg{\rangle} = \frac{1}{\kappa} \bigg{\langle} \bigg{|} \frac{1}{\kappa} - \frac{4l+2}{\sum_{j=1}^{2l+1}N_j}\bigg{|} \bigg{\rangle}$$ as $N \to \infty$. 
 Letting $C(l)$ be defined by $$C(l)=\bigg{\langle} \bigg{|} \frac{1}{\kappa} - \frac{4l+2}{\sum_{j=1}^{2l+1}N_j}\bigg{|} \bigg{\rangle},$$ we will now show that $\lim_{l \to \infty}C(l)=0$. 
 $N_j$ is a local function of the figures in the graph, and $\frac{1}{2l+1} \sum_{j=1}^{2l+1}N_j$ is an ergodic average of a local function of the graph so \textbf{(G1)} implies $\frac{1}{2l+1} \sum_{j=1}^{2l+1}N_j$ converges to $2\kappa$. Therefore for any $\varepsilon > 0$ we have that
$$\lim_{l \to \infty} \text{P} \biggl( \bigg{|}\frac{1}{\kappa} - \frac{4l+2}{\sum_{j=1}^{2l+1}N_j}\bigg{|} > \varepsilon \biggr) = 0.$$ Since $1 \leq \kappa \leq 2$ and $2 \leq N_j \leq 4$, the expression  $$ \bigg{|} \frac{1}{\kappa} - \frac{4l+2}{\sum_{j=1}^{2l+1}N_j} \bigg{|} $$ is bounded by a constant and we can apply the dominated convergence theorem to conclude that $\lim_{l \to \infty}C(l)=0$. Convergence of $\frac{1}{\kappa}$ to $\frac{4l+2}{\sum_{j=1}^{2l+1}N_j}$ in probability, instead of the stronger condition of almost sure convergence, is sufficient for the dominated convergence theorem to be applied \cite[p.258]{shiryaev}. Hence the error term which arises in replacing $\bar{\omega}_j^l$ with $2\kappa \omega_j^l$ is $\frac{1}{\gamma}\tilde{H}_N(t)$ plus $O(N)C(l)$.

Replacing $N_j$ with $2\kappa$ in the remaining part of the expression \eqref{eq:7} creates the error term \[ \begin{aligned} \int_{\Omega_N} \sum_{j=1}^{T_N} ( N_j - 2\kappa ) & \rho(t,j/T_N) \Phi'(\rho(t,j/T_N))F(t,j/T_N) f_t^N  \text{d}\nu^N \\ &
\leq \sum_{j=1}^{T_N} ( N_j - 2\kappa ) \bigg{\|} \rho(t,x) \Phi'(\rho(t,x))F(t,x) \bigg{\|}_\infty \end{aligned} \]
 which is controlled by noting that by the ergodicity condition \textbf{(G1)}, $$\bigg{|}  \sum_{j=1}^{T_N}(N_j - 2\kappa) \bigg{|} = | 2N - 2\kappa T_N | = o(N).$$ 
By performing the above replacements, we have obtained that
\[ \begin{aligned} \int_{\Omega_N} \frac{\partial_t \tilde{\psi}_t^N}{\tilde{\psi}_t^N} \ f_t^N \text{d}\nu^N & = \int_{\Omega_N} \sum_{j=1}^{T_N} 2\kappa \Bigl( \omega_j^l- \rho(t,j/T_N) \Bigr) \Phi'(\rho(t,j/T_N))F(t,j/T_N) f_t^N   \text{d}\nu^N \\ & \quad + o(N) + O(N)C(l). \end{aligned} \]

Next,  the factor  $\frac{N^2}{T_N^2}$ in front of each term in the last line of \eqref{eq:5} must be replaced by $\kappa^2$, creating the error  term 
\[ \begin{aligned}
\int_{\Omega_N} & \bigg{(}\kappa^2 - \frac{N^2}{T_N^2}\bigg{)}\sum_{j=1}^{T_N} (g(\omega_{x_j,1})+g(\omega_{x_j,-1}))\cfrac{\partial_x^2 \Phi(\rho(t,j/T_N))}{\Phi(\rho(t,j/T_N))} f_t^N  \text{d}\nu^N \\ & \leq
 \bigg{|}\kappa^2 - \frac{N^2}{T_N^2}\bigg{|} \bigg{\|} \cfrac{\partial_x^2 \Phi(\rho(t,x))}{\Phi(\rho(t,x))} \bigg{\|}_\infty \int_{\Omega_N} \sum_{j=1}^{T_N}  (g(\omega_{x_j,1})+g(\omega_{x_j,-1})) f_t^N  \text{d}\nu^N \\ & \leq
\bigg{|}\kappa^2 - \frac{N^2}{T_N^2}\bigg{|} \bigg{\|} \cfrac{\partial_x^2 \Phi(\rho(t,x))}{\Phi(\rho(t,x))} \bigg{\|}_\infty \int_{\Omega_N} \sum_{j=1}^{T_N}  B(\omega_{x_j,1}+\omega_{x_j,-1})) f_t^N  \text{d}\nu^N
\\ & \leq
\bigg{|}\kappa^2 - \frac{N^2}{T_N^2}\bigg{|} \bigg{\|} \cfrac{\partial_x^2 \Phi(\rho(t,x))}{\Phi(\rho(t,x))} \bigg{\|}_\infty \int_{\Omega_N} \sum_{j=1}^{T_N}  B \hat{\omega}_j f_t^N  \text{d}\nu^N
 \end{aligned} \]
where the condition \textbf{(SLG)}, which implies that $g(\omega_j) \leq B \omega_j$ for some $B > 0$, is used in the third line. This error term is $o(N)$ due to Lemma \ref{lm:2} and \textbf{(G1)}. Similarly, $\frac{N^2}{T_N^2}$ is replaced by $\kappa^2$ in \eqref{eq:6} creating the $o(N)$ error term \[ \begin{aligned}
 \bigg{(}2\kappa^2 & - \frac{2N^2}{T_N^2}\bigg{)}  \sum_{j=1}^{T_N}  \Phi(\rho(t,j/T_N)) \cfrac{\partial_x^2 \Phi(\rho(t,j/T_N))}{\Phi(\rho(t,j/T_N))}  \\ &
 \leq \sum_{j=1}^{T_N} \bigg{|}2\kappa^2 - \frac{2N^2}{T_N^2}\bigg{|} \bigg{\|} \Phi(\rho(t,x)) \cfrac{\partial_x^2 \Phi(\rho(t,x))}{\Phi(\rho(t,x))} \bigg{\|}_\infty. \end{aligned} \]

We now aim to replace $g(\omega_{x_j,1})+g(\omega_{x_j,-1})$ in the last line of \eqref{eq:5} with a more convenient expression. 
Define
\begin{equation*} g_j^l(\omega) = \frac{1}{4l+2}\sum_{|k-j| \leq l}(g(\omega_{x_{k},1})+g(\omega_{x_{k},-1})) \end{equation*}  the average sum of jump rates from the centres of the tiles near $j$. First $g(\omega_{x_j,1})+g(\omega_{x_j,-1})$ is replaced by $2 g_j^l(\omega)$, giving an error identical to \eqref{eq:9} except that the expression contains $g(\omega_{x_j,1})+g(\omega_{x_j,-1})$ instead of $\hat{\omega}_j$. $2g_j^l(\omega)$ is then replaced by $2\Phi(\omega_j^l)$, creating an error which is controlled using the following lemma.

\begin{lemma}[One block estimate]\label{lm:3}
 $$\limsup_{l \to \infty} \limsup_{N \to \infty} \frac{1}{N} \int_{\Omega_N} \sum_{j=1}^{T_N} |g_j^l(\omega)-\Phi(\omega_j^l)|\text{d}\mu_t^N = 0.$$
\end{lemma}
\begin{proof}
 This lemma is proved in Section \ref{sc:5}.
\end{proof}

Altogether, to obtain \eqref{eq:8},  \begin{equation*}
\begin{aligned} &
\int_{\Omega_N}  \frac{1}{\tilde{\psi}_t^N}(N^2 L_N^* \tilde{\psi}_t^N - \partial_t \tilde{\psi}_t^N )f_t^N \text{d}\nu^N 
\\ & = \int_{\Omega_N} \sum_{j=1}^{T_N} \Bigg{\{} \frac{N^2}{T_N^2} (g(\omega_{x_j,1})+g(\omega_{x_j,-1}))\frac{\partial_x^2 \Phi(\rho(t,j/T_N))}{\Phi(\rho(t,j/T_N))}\\
 &\quad- \frac{2N^2}{T_N^2} \Phi(\rho(t,j/T_N)) \frac{\partial_x^2 \Phi(\rho(t,j/T_N))}{\Phi(\rho(t,j/T_N))} 
\\ &\quad - \Bigl( \hat{\omega}_j-N_j\rho(t,j/T_N) \Bigr) \frac{\kappa \Phi'(\rho(t,j/T_N))\partial_x^2 \Phi(\rho(t,j/T_N))}{\Phi(\rho(t,j/T_N))} \Bigg{\}} f_t^N \text{d}\nu^N 
\\ & \leq
\int_{\Omega_N} \sum_{j=1}^{T_N} 2 \kappa F(t,j/T_N) \Big{\{} \Phi(\omega_j^l) - \Phi(\rho(t,j/T_N))\\
 &\quad- \Phi'(\rho(t,j/T_N))( \omega_j^l - \rho(t,j/T_N)) \Big{\}}f_t^N \text{d}\nu^N  
\\ & \quad +  \frac{2l(2l+1)}{N} \Big{\|} \partial_x(F(t,x)\Phi'(\rho(t,x))) \Big{\|}_\infty \int_{\Omega_N} \sum_{j=1}^{T_N} \hat{\omega}_j f_t^N \text{d}\nu^N 
\\ & \quad + \frac{1}{\gamma}\tilde{H}_N(t)\\
 &\quad+\frac{1}{\gamma} \sum_{j=1}^{T_N} \log \mathbb{E}_{\tilde{\nu}_{\rho(t,.)}^N} \bigg{[} \exp \Bigl( \frac{\| \kappa F(t,x)\Phi'(\rho(t,x)) \|_\infty \gamma \hat{\omega}_j}{2l+1} \sum_{|k - j| \leq l} \Big{|} \frac{1}{\kappa }-\frac{4l+2}{\sum_{|m-k| \leq l} N_m}  \Big{|} \Bigr) \bigg{]}
\\ & \quad +  \int_{\Omega_N} \sum_{j=1}^{T_N} ( N_j - 2\kappa )  \rho(t,j/T_N) \Phi'(\rho(t,j/T_N))F(t,j/T_N) f_t^N  \text{d}\nu^N
\\ & \quad + \int_{\Omega_N}  \bigg{(}\kappa^2 - \frac{N^2}{T_N^2}\bigg{)}\sum_{j=1}^{T_N} (g(\omega_{x_j,1})+g(\omega_{x_j,-1}))\cfrac{\partial_x^2 \Phi(\rho(t,j/T_N))}{\Phi(\rho(t,j/T_N))} f_t^N  \text{d}\nu^N
\\ & \quad +  \kappa^2 \int_{\Omega_N} \sum_{j=1}^{T_N} |g_j^l(\omega)-\Phi(\omega_j^l)| \Big{\|} \Phi'(\rho(t,j/T_N))F(t,j/T_N) \Big{\|}_\infty f_t^N \text{d}\nu^N
\\ & \quad +  \bigg{(}2\kappa^2  - \frac{2N^2}{T_N^2}\bigg{)}  \sum_{j=1}^{T_N}  \Phi(\rho(t,j/T_N)) \cfrac{\partial_x^2 \Phi(\rho(t,j/T_N))}{\Phi(\rho(t,j/T_N))}
\\ = & \int_{\Omega_N} \sum_{j=1}^{T_N} 2 \kappa F(t,j/T_N) \Big{\{} \Phi(\omega_j^l) - \Phi(\rho(t,j/T_N)) - \Phi'(\rho(t,j/T_N))( \omega_j^l - \rho(t,j/T_N)) \Big{\}}f_t^N \text{d}\nu^N
\\ & \quad + \frac{1}{\gamma}\tilde{H}_N(t) + o(N) + O(N)C(l)
\end{aligned} \end{equation*}
since we showed that the error terms in the middle equality consist of an $o(N)$ error and an $O(N)C(l)$ error.

 \subsection{Bounding $\tilde{H}_N(t)$}\label{sc:4.4}

Applying the entropy inequality to the main expression with $\gamma > 0$,
\begin{equation} \begin{aligned}
\int_{\Omega_N} & \sum_{j=1}^{T_N} 2 \kappa F(t,j/T_N)  \Big{\{} \Phi(\omega_j^l) - \Phi(\rho(t,j/T_N)) - \Phi'(\rho(t,j/T_N))( \omega_j^l - \rho(t,j/T_N)) \Big{\}}f_t^N \text{d}\nu^N \\ & \quad + \frac{1}{\gamma}\tilde{H}_N(t) \\ & 
\leq \frac{1}{\gamma} \log \int_{\Omega_N} \exp \Bigl(  \sum_{j=1}^{T_N} 2 \kappa \gamma F(t,j/T_N) \Big{\{} \Phi(\omega_j^l) - \Phi(\rho(t,j/T_N)) \\ & \quad - \omega_j^l \Phi'(\rho(t,j/T_N)) +  \rho(t,j/T_N) \Phi'(\rho(t,j/T_N))\Big{\}} \Bigr) \text{d}\tilde{\nu}_{\rho(s,.)}^N + \frac{2}{\gamma}\tilde{H}_N(t) + o(N).
\end{aligned}  \label{eq:gamma} \end{equation}
The coefficient of $\tilde{H}_N(t)$ is $2/\gamma$ instead of $1/\gamma$ since there are two separate applications of the entropy inequality.

Let
 \begin{equation}
G(u, \omega_j^l) = 2 \kappa \gamma F(t,u) \Big{\{} \Phi(\omega_j^l) - \Phi(\rho(t,u)) - (\omega_j^l-\rho(t,u)) \Phi'(\rho(t,u)) \Big{\}}. \label{eq:10}
\end{equation} 
 By Lipschitz continuity of $\Phi$ on $[0,\lambda]$ there exists $C_0$ for which $\Phi(\lambda) - \Phi(0) = \Phi(\lambda) \leq \lambda C_0$ and $\lambda \Phi(\rho(s,u)) \leq \lambda C_0$. Recall that $K_2$ is an upper bound for the value of $\rho(t,x)$ when $t = 0$, and hence for all $t > 0$ by the Maximum Principle. Defining
$$C_1 = 4\kappa \gamma \| F \|_\infty (\sup_{\beta \in [0,K_2]}\Phi(\beta)+C_0 K_2), \ C_2 = 8\kappa \gamma \| F\|_\infty C_0,$$ we have that $|G(u,\lambda)| \leq C_1 + C_2 \lambda$.

The first term on the right hand side of \eqref{eq:gamma} is equal to
\begin{equation*}
\frac{1}{\gamma}\log \int_{\Omega_N} \exp \Bigl( \sum_{j=1}^{T_N}G(j/T_N,\omega_j^l) \Bigr) \text{d}\tilde{\nu}_{\rho(s,.)}^N.
\end{equation*}
By repeated applications of H\"older's inequality, with $T_N = (2l+1)k$ for some integer $k$, \begin{equation}  \begin{aligned}
\log \int_{\Omega_N} & \exp \Bigl( \sum_{j=1}^{T_N} G(j/T_N,\omega_j^l)\Bigr) \text{d}\tilde{\nu}_{\rho(s,.)}^N \\
& = \log \int_{\Omega_N} \biggl( \prod_{i =0}^{2l} \prod_{m=0}^{k-1} \exp G\Bigl( \frac{(2l+1)m+i}{T_N},\omega_{(2l+1)k+i}^l \Bigr) \biggr) \text{d}\tilde{\nu}_{\rho(s,.)}^N
\\ & \leq \log \prod_{i =0}^{2l} \biggl( \int_{\Omega_N} \Bigl( \prod_{m=0}^{k-1} \exp \Bigl( G\Bigl( \frac{(2l+1)m+i}{T_N},\omega_{(2l+1)m+i}^l \Bigr) \Bigr)^{2l+1}\text{d}\tilde{\nu}_{\rho(s,.)}^N \biggr)^{\frac{1}{2l+1}} \\
\\ & = \log \prod_{i = 0}^{2l} \prod_{m=0}^{k-1}  \biggl( \int_{\Omega_N} \Bigl( \exp  \Bigl( G \Bigl( \frac{(2l+1)m+i}{T_N},\omega_{(2l+1)m+i}^l \Bigr) \Bigr) \Bigr)^{2l+1}\text{d}\tilde{\nu}_{\rho(s,.)}^N \biggr)^{\frac{1}{2l+1}}
\\ & = \frac{1}{2l+1}\sum_{j=1}^{T_N} \log \int_{\Omega_N} \exp \Bigl( (2l+1)G(j/T_N,\omega_j^l)\Bigr) \text{d}\tilde{\nu}_{\rho(s,.)}^N
\end{aligned} \label{eq:11} \end{equation}
where the independence of $\omega_j^l$ that are a distance $2l+1$ tiles apart is used in the last line. The condition that $T_N = (2l+1)k$ can be assumed without loss of generality since if $T_N = (2l+1)k+r$ for some $0 < r < 2l+1$ then H\"older's inequality can be applied to the first $(2l+1)k$ terms and the remaining $r$ can each be bounded by a constant. This is achieved by bounding $G$ as a linear function of the $\omega_j^l$ and letting $\gamma$ be small enough that the Laplace transform of each $\omega_{j, 1}$ and $\omega_{j,-1}$ is finite.

Every chain of $2l+1$ tiles contains between $4l+2$ and $8l+4$ vertices. It will later be necessary to decompose the sum based on the number of vertices in the $2l+1$ tiles surrounding each $j$, writing it as
$$\frac{1}{2l+1} \sum_{m=4l+2}^{8l+4} \sum_{j \in S_m} \log \int_{\Omega_N} \exp \Bigl(  (2l+1)G(j/T_N,\omega_j^l)\Bigr) \text{d}\tilde{\nu}_{\rho(s,.)}^N$$
where $$S_m= \Big{\{} j \in \{ 1, \cdots, T_N \} : \sum_{|k-j| \leq l}N_k = m \Big{\}}.$$

For a fixed $m \in \{ 4l+2, \cdots, 8l+4 \}$ let $j_{m} \in \mathbb{T}_N$ be the site of the centre closest to the centre of tile $j$ such that the $2l+1$ tiles surrounding $j_{m}$ have $m$ vertices. 
In the case where there are two tile centres equally close to $j$ such that both are surrounded by $2l+1$ tiles with $m$ vertices  we let $j_m$ be  the one to the right of $j$.
For sufficiently large $N$ this is well defined since sequences of $2l+1$ tiles of all possible sizes occur in the graph. For any $u \in \mathbb{T}$ the marginals $\tilde{\nu}_{\rho(s,\lfloor u T_N \rfloor_m / T_N)}^1$ of $\tilde{\nu}_{\rho(s,.)}^N$ converge to the probability measure $\mathbb{P}_{\rho(s,u)}$ defined in \eqref{eq:4}. This is due to the continuity of $\rho$ and the fact that by \textbf{(G1)}, for all $\gamma \in [0,1]$, $\frac{ \lfloor \gamma T_N \rfloor_m}{T_N} \to \gamma$ as $N \to \infty$. Since $l$ is fixed, $\exp((2l+1)  G(\lfloor u T_N \rfloor_m / T_N, \omega_{ \lfloor u T_N \rfloor_m}^l))$ is a bounded cylinder function. Hence,
for all $u \in \mathbb{T}$ and $m \in \{ 4l+2, \cdots, 8l+4\}$, 
\begin{equation} \int_{\Omega_N} \exp \Bigl(  (2l+1) G(\lfloor u T_N \rfloor_m / T_N, \omega_{\lfloor u T_N \rfloor_m}^l \Bigr) \text{d}\nu_{\rho(s,u)}^N \to \mathbb{E}_{\rho(s,u)}\Big{[} \exp \Bigl( (2l+1) G \Bigl( u,\frac{1}{m}\sum_{j=1}^{m}X_j \Bigr) \Bigr) \Big{]} \label{eq:12} \end{equation}
as $N \to \infty$, where the $X_j$ are independent, identically distributed random variables with distribution $\mathbb{P}_{\rho(s,u)}$ as in \eqref{eq:4}. 
The decomposition of the sum into separate values of $m$ is necessary here since the random variables $\omega_{\lfloor u T_N \rfloor}^l$ contain varying numbers of terms and converge to different limiting random variables depending on this number.

For a fixed $l$, we now replace the integral with respect to $\tilde{\nu}_{\rho(s,.)}^N$ in the last line of \eqref{eq:11} with a sum of integrals over $\mathbb{T}$ for each $m \in \{ 4l+2, \cdots, 8l+4 \}$. This is achieved by first replacing it with an integral with respect to $\nu_{\rho(s,.)}^N$ and then applying the above decomposition. 
We have that
\begin{equation*} \begin{aligned}
 \frac{1}{ (2l+1)N} & \sum_{j=1}^{T_N}  \log \mathbb{E}_{\tilde{\nu}_{\rho(s,.)}^N} \Big{[} \exp ( (2l+1) G( j/T_N,\omega_j^l) \Big{]} \\
& \leq  \frac{1}{(2l+1)\kappa } \sum_{m=4l+2}^{8l+4} p(l,m) \int_{\mathbb{T}} \log \mathbb{E}_{\rho(s,x)} \Big{[} \exp \Bigl(  (2l+1) G \Bigl( x,\frac{1}{m}\sum_{j=1}^{m}X_j \Bigr) \Bigr) \Big{]} \text{d}x \\
& \quad + \Bigg{|} \frac{1}{ (2l+1)N}\sum_{j=1}^{T_N} \log \mathbb{E}_{\tilde{\nu}_{\rho(s,.)}^N} \Big{[} \exp ( (2l+1) G( j/T_N,\omega_j^l) \Big{]} \\
& \qquad - \frac{1}{ (2l+1)N}\sum_{j=1}^{T_N} \log \mathbb{E}_{\nu_{\rho(s,.)}^N}  \Big{[} \exp ( (2l+1) G(x_j/N,\omega_j^l)\Big{]} \Bigg{|} \\
& \quad  + \Bigg{|} \frac{1}{ (2l+1)N} \sum_{m=4l+2}^{8l+4} \sum_{j \in S_m} \log \mathbb{E}_{\nu_{\rho(s,.)}^N}  \Big{[} \exp ( (2l+1) G(x_j/N,\omega_j^l))\Big{]}   \\
& \qquad - \frac{1}{ (2l+1)\kappa } \sum_{m=4l+2}^{8l+4} p(l,m) \int_{\mathbb{T}} \log \mathbb{E}_{\rho(s,x)} \Big{[} \exp \Bigl(  (2l+1) G\Bigl( x,\frac{1}{m}\sum_{j=1}^{m}X_j \Bigr) \Bigr) \Big{]} \text{d}x \Bigg{|}  
\end{aligned} \end{equation*}
where $p(l,m)$ is the limiting proportion of sequences of $2l+1$ tiles with $m$ vertices as $N \to \infty$. For a fixed $l$, $\sum_{m=4l+2}^{8l+4}p(l,m) = 1$ and $\lim_{N \to \infty}\frac{|S_m|}{T_N} = p(l,m)$. 
The integrals over $\mathbb{T}$ are multiplied by $\frac{p(l,m)}{\kappa}$ since each Riemann sum contains $|S_m|$ terms but is divided by $N$.

The first error term is dealt with using a coupling of the measures $\tilde{\nu}_{\rho(s,.)}^N$ and $\nu_{\rho(s,.)}^N$. 
Let $f(\omega) = (2l+1) |G( \lfloor jT_N \rfloor/T_N,\omega_j^l)| \leq (2l+1)[C_1 + \omega_j^l C_2].$ 
Define a coupling measure $\mathbb{P}$ for two random variables $\omega$ and $\omega'$ with the property that $\omega_x \leq \omega_x'$ at each vertex $x$ with probability 1, with $\omega_x$ distributed according to the minimum of the densities of the two marginals of $\nu_{\rho(t,.)}^N$ and  $\tilde{\nu}_{\rho(t,.)}^N$ at vertex $x$ and $\omega_x'$ according to the maximum. This coupling measure exists due to Lemma A.2 in \cite{balazs} and because $g$ is increasing. Let $\mathbb{E}$ be the expectation with respect to $\mathbb{P}$.
We have that
\[ \begin{aligned}
\mathbb{E}[e^{f(\omega')}-e^{f(\omega)}]  & = \mathbb{E} \Big{[} \Bigl( e^{f(\omega')}-e^{f(\omega)}\Bigr)\mathbf{1}\Big{\{} \cup_{|k - j | \leq l} \cup_{s \in \text{tile} \ k} \omega_s' \neq \omega_s \Big{\}} \Big{]} \\ & \leq
\sum_{|k - j| \leq l} \sum_{s \in \text{tile} \ k} \mathbb{E}\Big{[} \Bigl( e^{f(\omega')}-e^{f(\omega)} \Bigr) \mathbf{1}\{  \omega_s' \neq \omega_s \} \Big{]} \\ &
\leq \sum_{|k - j| \leq l} \sum_{s \in \text{tile} \ k} \mathbb{E} \Big{[} \Bigl( e^{f(\omega')}+e^{f(\omega)}\Bigr)^p \Big{]}^{1/p} \mathbb{P}( \omega_s' \neq \omega_s )^{1/q} \\ &
\leq \sum_{|k-j| \leq l} \sum_{s \in \text{tile} \ k}
\Bigl( 2^{p-1}\mathbb{E}[e^{pf(\omega')} + e^{pf(\omega)}]\Bigr)^{1/p}(\mathbb{P}( \omega_s' \neq \omega_s ))^{1/q} \\ &
\leq \sum_{|k - j| \leq l} \sum_{s \in \text{tile} \ k} C(\rho)\mathbb{E}[\omega_s'-\omega_s]^{1/q}
\end{aligned}
\]
for some exponents $p$ and $q$ such that $\frac{1}{p} + \frac{1}{q} = 1$ and $p$ is small enough that $\mathbb{E}[e^{pf(\omega)}]$ converges. Sending $N \to \infty$ while $l$ is fixed, the difference between the parameters given by $\nu_{\rho(t,.)}$ and  $\tilde{\nu}_{\rho(t,.)}$ converges to 0 by \textbf{(G1)}. Thus $\mathbb{E}[\omega_s'-\omega_s]$ tends to 0 for every vertex $s$, and there are at most $8l+4$ vertices for a fixed value of $l$ so the expression converges to 0 as $N \to \infty$.

This coupling can be applied to the error term since the expressions $\mathbb{E}_{\tilde{\nu}_{\rho(s,.)}^N}[\exp((2l+1)G(j/T_N,\omega_j^l)]$ and $\mathbb{E}_{\nu_{\rho(s,.)}^N}[\exp((2l+1)G(x_j/N,\omega_j^l)]$ are bounded away from 0 and thus $\log$ is continuous. The bound depends on $j$ only through the parameter values from $\nu_{\rho(t,.)}^N$ and  $\tilde{\nu}_{\rho(t,.)}^N$ and the difference between these goes to 0 uniformly.

The second error term vanishes due to the convergence of the sum to a Riemann integral, using the fact that $\lim_{N \to \infty}\frac{|S_m|}{T_N} = p(l,m)$ and the pointwise convergence given by \eqref{eq:12}.

Now, \begin{equation*}  \begin{aligned} & \frac{1}{(2l+1)\kappa } \sum_{m=4l+2}^{8l+4} p(l,m) \int_{\mathbb{T}} \log \mathbb{E}_{\rho(s,x)} \Big{[} \exp \Bigl(  (2l+1) G \Bigl( x,\frac{1}{m}\sum_{j=1}^{m}X_j \Bigr) \Bigr) \Big{]} \text{d}x  \\ &  \leq
\frac{1}{(2l+1)\kappa} \sup_{m \in \{4l+2, \cdots, 8l+4 \}} \int_{\mathbb{T}} \log \mathbb{E}_{\rho(s,x)} \Big{[} \exp \Bigl(  (2l+1) G \Bigl( x,\frac{1}{m}\sum_{j=1}^{m}X_j \Bigr) \Bigr) \Big{]} \text{d}x \end{aligned} \end{equation*} so from now on we will consider $$\frac{1}{(2l+1)\kappa} \int_{\mathbb{T}} \log \mathbb{E}_{\rho(s,x)} \Big{[} \exp \Bigl(  (2l+1) G \Bigl( x,\frac{1}{m}\sum_{j=1}^{m}X_j \Bigr) \Bigr) \Big{]} \text{d}x$$ for some $m \in \{ 4l+2, \cdots, 8l+4 \}$.

\begin{proposition}\label{pr:4}
\phantom{a}
\begin{enumerate} 
 \item[\textbf{(i)}]\begin{equation*} \begin{aligned} & 
\limsup_{l \to \infty} \frac{1}{(2l+1)\kappa} \int_{\mathbb{T}}   \log \mathbb{E}_{\rho(s,u)} \Big{[} \exp \Bigl(  (2l+1) G \Bigl( u,\frac{1}{m}\sum_{j=1}^{m}X_j \Bigr) \Bigr) \Big{]} \text{d}u \\ & \quad
\leq \int_{\mathbb{T}}\sup_{\lambda > 0}((\gamma | F(s,u)M(\lambda,\rho(s,u)) | -J_{\rho(s,u)}(\lambda))\text{d}u
\end{aligned}
\end{equation*}
where $J_{\rho(s,u)}$ is the large deviations rate function corresponding to iid random variables distributed according to $\mathbb{P}_{\rho(s,u)}$ and $M(\lambda,\rho)=2\kappa \{ \Phi(\lambda)-\Phi(\rho)-(\lambda-\rho) \Phi'(\rho) \}$.
\item[\textbf{(ii)}] There exists $\gamma > 0$ such that for all $u \in \mathbb{T}$, $$\sup_{\lambda > 0}(\gamma | F(s,u)M(\lambda, \rho(s,u))| -J_{\rho(s,u)}(\lambda)) \leq 0.$$
\end{enumerate}
\end{proposition}

\begin{proof}
 First we prove \textbf{(i)}.
A cutoff function is used, since when taking the limit as $l \to \infty$ $G$ may be unbounded. For $A > 0$ let
$$G_A(u,\lambda)=G(u,\lambda)\mathbf{1}\{ |\lambda | \leq A \} + G(u,A)\mathbf{1}\{ | \lambda | > A \}.$$
Then,
\begin{equation} \begin{aligned} &
\frac{1}{(2l+1)\kappa} \int_{\mathbb{T}} \log \mathbb{E}_{\rho(s,u)}\Big{[} \exp \Bigl( (2l+1)G\Bigl( u,\frac{1}{ m }\sum_{j=1}^{m}X_j  \Bigr) \Big{]} \ \text{d}u \\  & 
\leq \frac{1}{(2l+1)\kappa} \int_{\mathbb{T}} \log \bigg{\{} \mathbb{E}_{\rho(s,u)}\Big{[} \exp \Bigl( (2l+1)G_A \Bigl( u,\frac{1}{m}\sum_{j=1}^{m}X_j \Bigr) \Bigr)  \Big{]} \\ & \quad
+ e^{C_1 l} \mathbb{E}_{\rho(s,u)}\Big{[} \mathbf{1}\Big{\{}\frac{1}{ m}\sum_{j=1}^{m}X_j \geq A \Big{\}}\Big{]}^{1/q} \mathbb{E}_{\rho(s,u)}\Big{[} \exp \Bigl( pC_2 \sum_{j=1}^{m} X_j \Bigr) \Big{]}^{1/p} \bigg{\}} \ \text{d}u 
\end{aligned} \label{eq:13}
\end{equation}
where the second term is an upper bound for the difference between the $G$ and $G_A$ terms, to which H\"older's inequality is applied with $\frac{1}{p}+\frac{1}{q}=1$. Considering the first term, by monotonicity,
\begin{equation*} \begin{aligned} & \frac{1}{2l+1}  \int_{\mathbb{T}} \log  \mathbb{E}_{\rho(s,u)}\bigg{[} \exp \Bigl( (2l+1)G_A\Bigl( u,\frac{1}{m}\sum_{j=1}^{m}X_j \Bigr) \Bigr)  \bigg{]} \ \text{d}u \\ & \leq \frac{1}{2l+1} \int_{\mathbb{T}} \log \mathbb{E}_{\rho(s,u)}\bigg{[} \exp \Bigl( m \Big{|} G_A\Bigl( u,\frac{1}{m}\sum_{j=1}^{m}X_j \Bigr) \Bigr) \Big{|} \bigg{]} \ \text{d}u. \end{aligned} \end{equation*}
Since $G_A$ is bounded and $4l + 2 \leq m$, as $l \to \infty$, $m \to \infty$ and this term converges to
$$ \int_{\mathbb{T}} \sup_{\lambda > 0}(|G_A(u,\lambda)|-J_{\rho(s,u)}(\lambda)) \ \text{d}u $$
 by the Laplace-Varadhan Lemma.

To show that the second term inside the expectation in \eqref{eq:13} vanishes as $A \to \infty$, note that 
$$\frac{1}{2l+1}\log \mathbb{E}_{\rho(s,u)}\Big{[} \exp \Bigl( pC_2 \sum_{j=1}^{m}X_j \Bigr) \Big{]} = \frac{m}{2l+1} \log \mathbb{E}_{\nu_{\rho(s,u)}^1}\Big{[} \exp (pC_2 \omega_0 ) \Big{]}$$ by independence, and
$$\mathbb{E}_{\nu_{\rho(s,u)}^1} \Big{[} \exp(p C_2 \omega_0 ) \Big{]} = \sum_{k=0}^\infty \frac{(e^{pC_2})^k \Phi(\rho(s,u))^k}{g(k)! Z(\Phi(\rho(s,u)))} = \frac{Z(e^{pC_2}\Phi(\rho(s,u)))}{Z(\Phi(\rho(s,u)))}.$$
Let $\phi^*$ be a constant less than or equal to the radius of convergence of the power series $Z(\phi)$ and such that $\phi^* > \Phi(K_2)$, and recall that $C_2$ was defined as $C_2 = 8\kappa \gamma \|F\|_\infty C_0$. 
If $\gamma < \cfrac{\log(\phi^*/\Phi(K_2))}{8\kappa \|F\|_\infty C_0}$ then we have $C_2 < \log \Bigl( \frac{\phi^*}{\Phi(K_2)}\Bigr)$ and
$$Z(e^{pC_2}\Phi(\rho(s,u))) < Z(e^{p\log (\phi^*/ \Phi(K_2))}\Phi(\rho(s,u))) = Z \Bigl( \Bigl( \frac{\phi^*}{\Phi(K_2)} \Bigr)^p \Phi(\rho(s,u))\Bigr).$$  The term inside $Z$ is less than the radius of convergence $\phi^*$ if $p$ is chosen such that $$p < \frac{\log \phi^* - \log \Phi(\rho(s,u))}{\log \phi^* - \log \Phi(K_2)}.$$
A large deviations principle can be used to show that the factor $\mathbb{P}_{\rho(s,u)}(\frac{1}{m}\sum_{i=1}^{m}X_j \geq A)$ converges to zero after sending $A \to \infty$.
Altogether, the integral of the right hand side of \eqref{eq:13}, divided by $2l+1$ and then sending $l \to \infty$, is now bounded by
$$\sup_{\lambda > 0}(|G_A(u,\lambda)|-J_{\rho(s,u)}(\lambda)),$$ an expression which is equal to $$\sup_{\lambda > 0}(|G(u,\lambda)|-J_{\rho(s,u)}(\lambda))$$ after sending $A \to \infty$ by the argument given in \cite[p.126-127]{landim}. 
We have established that \[ \begin{aligned} \limsup_{l \to \infty} \frac{1}{(2l+1)\kappa} & \int_{\mathbb{T}} \log \mathbb{E}_{\rho(s,u)}\Big{[} \exp \Bigl( (2l+1)G\Bigl( u,\frac{1}{ m }\sum_{j=1}^{m}X_j  \Bigr) \Big{]} \ \text{d}u
\\  & \leq \int_{\mathbb{T}} \sup_{\lambda > 0}(|G(u,\lambda)|-J_{\rho(s,u)}(\lambda)) \ \text{d}u
 \\ &
= \int_{\mathbb{T}} \sup_{\lambda > 0}(\gamma |F(s,u)M(\lambda,\rho(s,u))|-J_{\rho(s,u)}(\lambda)) \ \text{d}u. \end{aligned} \]
 This proves Proposition \ref{pr:4}(i). 
By identical arguments to those in Kipnis and Landim \cite{landim}, a constant $\gamma$ can be found such that $$\sup_{\lambda > 0}(\gamma |F(s,u)M(\lambda, \rho(s,u))|-J_{\rho(s,u)}(\lambda)) \leq 0,$$ which proves Proposition \ref{pr:4}(ii). 
\end{proof}

\begin{remark}
 We will need the following lemma in order to bound $\tilde{H}_N(t)$, which is analogous to a calculation in Funaki, Uchiyama and Yau \cite[p.8]{funaki}. The theorem is that if the initial relative entropy of the system at $t=0$ with respect to the normal reference measure $\nu_{\rho(t,.)}^N$ is small then it remains small with respect to the perturbed measure $\tilde{\nu}_{\rho(t,.)}^N$.
\end{remark}

\begin{lemma}\label{lm:5}
 If $H_N(0) = o(N)$ then $\tilde{H}_N(0)=o(N)$.
\end{lemma}

\begin{proof}
Since $T_N \leq N$, $\frac{1}{N}|\tilde{H}_N(0) -H_N(0)| \leq \frac{1}{T_N}|\tilde{H}_N(0) -H_N(0)|$. We have that
 \begin{equation*} 
\begin{aligned} \frac{1}{T_N} & |\tilde{H}_N(0) -H_N(0)| \\ &  = \frac{1}{T_N} \bigg{|} \int_{\Omega_N} f_0^N(\omega)\log  \psi_0^N(\omega) \ \text{d}\nu^N - \int_{\Omega_N} f_0^N (\omega)\log \tilde{\psi}_0^N(\omega) \ \text{d}\nu^N \bigg{|}
\\ & = \frac{1}{T_N}\bigg{|} \int_{\Omega_N} f_0^N(\omega)\log \Bigl( \frac{\psi_0^N(\omega)}{\tilde{\psi}_0^N(\omega)}\Bigr) \ \text{d}\nu^N \bigg{|}
\\  & =  \frac{1}{T_N} \Bigg{|} \int_{\Omega_N}f_0^N (\omega)\log \biggl( \prod_{j = 1}^{T_N} \cfrac{\Phi(\rho(0,\frac{x_j}{N}))^{\omega_{x_j,1}+\omega_{x_j,-1}}Z(\Phi(\rho(0,\frac{j}{T_N})))^2}{\Phi(\rho(0,\frac{j}{T_N})^{\omega_{x_j,1}+\omega_{x_j,-1}}Z(\Phi(\rho(0,\frac{x_j}{N})))^2} \biggr) \\ & \qquad
  + f_0^N(\omega)\log \biggl( \prod_{y \in \mathbb{T}_N \backslash \{x_j\}}  \cfrac{\Phi(\rho(0,\frac{y}{N}))^{\omega_{y,1}+\omega_{y,-1}}Z(\Phi(\rho(0,\frac{k_{y,1}}{T_N})))Z(\Phi(\rho(0,\frac{k_{y,-1}}{T_N}))))}{\Phi(\rho(0,\frac{k_{y,1}}{T_N})^{\omega_{y,1}}\Phi(\rho(0,\frac{k_{y,-1}}{T_N}))^{\omega_{y,-1}}Z(\Phi(\rho(0,\frac{y}{N})))^2}  \biggr) \text{d}\nu^N \Bigg{|} 
  \\ & = \frac{1}{T_N} \Bigg{|} \int_{\Omega_N} f_0^N(\omega) \sum_{j=1}^{T_N} \bigg{\{} (\omega_{x_j,1}+\omega_{x_j,-1})\log \biggl( \frac{\Phi(\rho(0,\frac{x_j}{N}))}{\Phi(\rho(0,\frac{j}{T_N}))} \biggr) + 2\log \biggl( \frac{Z(\Phi(\rho(0,\frac{j}{T_N})))}{Z(\Phi(\rho(0,\frac{x_j}{N})))} \biggr) \bigg{\}}
  \\ & \qquad + f_0^N(\omega)\sum_{y \in \mathbb{T}_N \backslash \{x_j\}} \bigg{\{} \omega_{y,1}\log \Bigl( \frac{\Phi(\rho(0,\frac{y}{N}))}{\Phi(\rho(0,\frac{k_{y,1}}{T_N}))} \Bigr)
  + \omega_{y,-1}\log \biggl( \frac{\Phi(\rho(0,\frac{y}{N}))}{\Phi(\rho(0,\frac{k_{y,-1}}{T_N}))} \biggr)
  \\ & \qquad + \log \biggl( \frac{Z(\Phi(\rho(0,\frac{k_{y,1}}{T_N})))}{Z(\Phi(\rho(0,\frac{y}{N})))} \biggr) + \log \biggl( \frac{Z(\Phi(\rho(0,\frac{k_{y,-1}}{T_N})))}{Z(\Phi(\rho(0,\frac{y}{N})))} \biggr) \bigg{\}} \text{d}\nu^N \Bigg{|}
  \end{aligned}
  \end{equation*}
 where $k_s$ is the tile number of vertex $s$. 
Now for each $j \in \{ 1, \cdots, T_N \}$, consider 
\begin{equation*} \begin{aligned} &  (\omega_{x_j,1}+\omega_{x_j,-1}) \log \Biggl( \frac{\Phi(\rho(0,\frac{x_j}{N}))}{\Phi(\rho(0,\frac{j}{T_N}))} \Biggr) 
\\ & = (\omega_{x_j,1}+\omega_{x_j,-1}) \log \Biggl( 1 +\cfrac{\Phi(\rho(0,\frac{x_j}{N})) - \Phi(\rho(0,\frac{j}{T_N}))}{\Phi(\rho(0,\frac{j}{T_N}))} \Biggr) 
\\ & = (\omega_{x_j,1}+\omega_{x_j,-1}) \Bigg{\{} \cfrac{\Phi(\rho(0,\frac{x_j}{N})) - \Phi(\rho(0,\frac{j}{T_N}))}{\Phi(\rho(0,\frac{j}{T_N}))} - \frac{1}{2} \Biggl( \cfrac{\Phi(\rho(0,\frac{x_j}{N})) - \Phi(\rho(0,\frac{j}{T_N}))}{\Phi(\rho(0,\frac{j}{T_N}))} \Biggr)^2 + o(1) \Bigg{\}}.
\end{aligned}
\end{equation*}
There are $T_N$ such terms and we can view their sum divided by $T_N$ as a Riemann sum and apply the dominated convergence theorem to show that the integrals converge to zero. Since $\rho$ is bounded above and below, the values taken by the step function are bounded by a constant. To establish piecewise convergence of the step functions to zero, observe that for any $\gamma \in \mathbb{T}$ and sequence of $j$ values such that $j/T_N \to \gamma$, the ergodicity condition \textbf{(G1)} implies that $|\frac{j}{T_N} - \frac{x_j}{N}| \to 0$ and therefore $|\Phi(\rho(0,\frac{j}{T_N})) - \Phi(\rho(0,\frac{x_j}{N}))| \to 0$ as $N \to \infty$.
The sum of the other terms can be shown to converge to 0 by the same argument.
\end{proof}

\begin{remark}
 Funaki, Uchiyama and Yau \cite{funaki} and Komoriya \cite{komoriya} also require a bound on their versions of $|\tilde{H}_N(0)-H_N(0)|$. However the summands in their expression are $O(\frac{1}{N})$, which is smaller than the bound we have in Lemma \ref{lm:5}.
\end{remark}

\begin{proposition}\label{pr:6}
 If $H_N(0) = o(N)$ then $\tilde{H}_N(t) = o(N)$ for any $t > 0$.
\end{proposition}

\begin{proof}
By Lemma \ref{lm:5} we have that $\tilde{H}_N(0) = o(N)$ and
\begin{equation*} \begin{aligned}
\tilde{H}_N(t) & = \int_0^t \partial_s \tilde{H}_N(s) \text{d}s + \tilde{H}_N(0) =  \int_0^t \partial_s \tilde{H}_N(s) \text{d}s + H_N(0)  + o(N) \\ &  \leq \int_0^t \int_{\Omega_N} \frac{1}{\tilde{\psi}_t^N} (N^2 L_N^* \tilde{\psi}_t^N - \partial_t \tilde{\psi}_t^N ) f_t^N \text{d} \nu^N \text{d}s + o(N) + O(N)C(l) \\ & 
=  \int_0^t \int_{\Omega_N} \sum_{j=1}^{T_N} F(s,j/T_N)M(\omega_j^l,j/T_N) \text{d}\mu_s^N \text{d}s + o(N) + O(N)C(l) \\ & 
\leq \frac{2}{\gamma} \int_0^t \tilde{H}_N(s)\text{d}s + \frac{1}{\gamma} \int_0^t \log \int_{\Omega_N} \exp \Bigl( \gamma \sum_{j=1}^{T_N} F(s,j/T_N)M(\omega_j^l,\rho(s,j/T_N)) \Bigr) \text{d}\tilde{\nu}_{\rho(s,.)}^N \text{d}s \\ & \qquad + o(N) + O(N)C(l). \end{aligned}
\end{equation*}
The third line is due to \eqref{eq:8} and in the fourth line the entropy inequality is applied and an additional  $\frac{1}{\gamma} \int_0^t \tilde{H}_N(s)\text{d}s$ term comes from the application of entropy inequality to another error term. 
Recalling the definition of $G$ given in \eqref{eq:10} and applying the repeated applications of H\"older's inequality in \eqref{eq:11} and finally the convergence to the Riemann integral, 
\[ \begin{aligned} & \frac{1}{\gamma N} \int_0^t \log \int_{\Omega_N} \exp \Bigl( \gamma \sum_{j=1}^{T_N} F(s,j/T_N)M(\omega_j^l,\rho(s,j/T_N)) \Bigr)\text{d}\tilde{\nu}_{\rho(s,.)}^N \text{d}s \\ & 
 \leq \frac{1}{(2l+1)N} \int_0^t  \sum_{j=1}^{T_N} \log \int_{\Omega_N} \exp \Bigl( (2l+1)G(j/T_N,\omega_j^l)\Bigr) \text{d}\tilde{\nu}_{\rho(s,.)}^N  \ \text{d}s + o(N) \\ & 
\leq \frac{1}{(2l+1)\kappa \gamma} \int_0^t \int_{\mathbb{T}}\log \mathbb{E}_{\rho(s,u)} \Big{[} \exp \Bigl( (2l+1)G\Bigl( u, \frac{1}{m}\sum_{j=1}^m X_j \Bigr) \Bigr) \Big{]} \text{d}u \ \text{d}s + o(N). \end{aligned}
\]
By Proposition \ref{pr:4}, the expression on the right hand side is $o(N)$ for sufficiently small $\gamma$ and we have that $$\tilde{H}_N(t) \leq \frac{2}{\gamma} \int_0^t \tilde{H}_N(s)\text{d}s + o(N) + O(N)C(l),$$
$\tilde{H}_N(t) \leq (o(N)+O(N)C(l))\exp(2t/\gamma)$ by Gronwall's inequality, $\tilde{H}_N(t) = o(N)+O(N)C(l)$. After sending $l \to \infty$, $\tilde{H}_N(t)=o(N)$.
\end{proof}

\subsection{Proof of Theorem \ref{tm:1}}\label{sc:4.5}

It remains to show the convergence in probability of the particle density to the solution of the partial differential equation. It is similar to Corollary 6.1.3 in Kipnis and Landim \cite[p.117-118]{landim}, with additional steps as a consequence of there being two reference measures. In \cite{funaki}, the analogous result is instead obtained via a large deviations estimate in Corollary 3.1.

\begin{proof}[Proof of Theorem \ref{tm:1}]
Let $\phi \in C^\infty(\mathbb{T})$ and let $\psi : \mathbb{Z}^{\{-1,1\} \times \mathbb{Z}} \to \mathbb{R}$ be a bounded cylinder function of $\omega$. Then convergence in probability can be shown by proving that
$$\lim_{N \to \infty}\mathbb{E}_{\mu_t^N}\Big{[}\Big{|}\frac{1}{N}\sum_{j \in \mathbb{T}_N} \phi(j/N)\psi(\tau_j \omega) - \int_{\mathbb{T}}\phi(u)\mathbb{E}_{\nu_{\rho(t,u)}}(\psi(\omega))\text{d}u \Big{|} \Big{]}= 0.$$ 
\textbf{Step 1.} 
By the triangle inequality,
\[ \begin{aligned}
\mathbb{E}_{\mu_t^N} \Big{|}\frac{1}{N} & \sum_{j \in \mathbb{T}_N} \phi(j/N)\psi(\tau_j \omega) - \int_{\mathbb{T}}\phi(u)\mathbb{E}_{\nu_{\rho(t,u)}}(\psi(\omega))\text{d}u \Big{|} \\ & \leq  \mathbb{E}_{\mu_t^N} \Big{|} \frac{1}{N}\sum_{j \in \mathbb{T}_N} \phi(j/N) \Bigl( \psi(\tau_j \omega) - \mathbb{E}_{\nu_{\rho(t,j/N)}}(\psi(\omega))\Bigr) \Big{|} \\ & \quad + \Big{|} \frac{1}{N}\sum_{j \in \mathbb{T}_N}\phi(j/N)\mathbb{E}_{\nu_{\rho(t,j/N)}}(\psi(\omega)) - \int_{\mathbb{T}}\phi(u)\mathbb{E}_{\nu_{\rho(t,u)}}(\psi(\omega))\text{d}u \Big{|}.
\end{aligned} \]
The second term converges to zero as it is a Riemann sum.

\noindent
\textbf{Step 2.} Considering the first term of the right hand side above and performing a one block estimate, 
\[ \begin{aligned}
&  \mathbb{E}_{\mu_t^N}  \Big{|} \frac{1}{N}\sum_{j \in \mathbb{T}_N} \phi(j/N) \Bigl( \psi(\tau_j \omega) - \mathbb{E}_{\nu_{\rho(t,j/N)}}(\psi(\omega))\Bigr) \Big{|}  \\  &
\leq \frac{1}{N}\sum_{j \in \mathbb{T}_N} |\phi(j/N)| \mathbb{E}_{\mu_t^N} \Big{|} \frac{1}{2l+1}\sum_{k : |k-j|\leq l} \Bigl( \psi(\tau_k \omega) - \mathbb{E}_{\nu_{\rho(t,k/N)}}(\psi(\omega))\Bigr) \Big{|} + 2\| \psi \|_{\infty}V(\phi,l/N) \\ &
\leq \frac{\| \phi \|_{\infty}}{N}\mathbb{E}_{\mu_t^N} \biggl( \sum_{j \in \mathbb{T}_N} \Big{|}\frac{1}{2l+1}\sum_{k:|k-j|\leq l}\Bigl( \psi(\tau_k \omega) - \mathbb{E}_{\nu_{\rho(t,j/N)}}(\psi(\omega)) \Bigr) \Big{|} \biggr) + o(1)
\end{aligned} \]
where $V(\phi,s)= \sup_{|x-y|\leq s}|\phi(x)-\phi(y)|$. Here, $ \| \psi \|_{\infty}V(\phi,l/N) = o(1)$ since $\psi$ is bounded and $\phi$ is continuous. 

\noindent
\textbf{Step 3.} 
Applying the entropy inequality, where $\gamma > 0$,
\begin{equation*} \begin{aligned} & \frac{1}{N}\mathbb{E}_{\mu_t^N}\Bigl( \sum_{j \in \mathbb{T}_N} \Big{|}\frac{1}{2l+1}\sum_{k:|k-j|\leq l}\Bigl( \psi(\tau_k \omega) - \mathbb{E}_{\nu_{\rho(t,j/N)}}(\psi(\omega)) \Bigr) \Big{|} \Bigr) \\ & \leq \frac{1}{\gamma N}\tilde{H}_N(t) + \frac{1}{\gamma N} \log \mathbb{E}_{\tilde{\nu}_{\rho(t,.)}^N} \exp \Bigl( \gamma \sum_{j \in \mathbb{T}_N} \Big{|} \frac{1}{2l+1} \sum_{k:|k-j|\leq l}\Bigl( \psi(\tau_k \omega)-\mathbb{E}_{\nu_{\rho(t,k/N)}}(\psi(\omega))\Bigr) \Big{|}\Bigr). \end{aligned}
\end{equation*}
\textbf{Step 4.} Let $d$ be the range of the cylinder function $\psi$. We assume $N$ is a multiple of $2l+d+1$, with $N = (2l+d+1)s$ for some integer $s$, since otherwise $N = (2l+d+1)s + r$ for some $0 < r < 2l + d + 1$ and H\"older's inequality can be applied to the first $(2l+d+1)s$ terms while the remaining $r$ terms are bounded by a constant. Applying H\"older's inequality to the second term $2l+d+1$ times,
\[ \begin{aligned} 
& \frac{1}{\gamma N}  \log \mathbb{E}_{\tilde{\nu}_{\rho(t,.)}^N} \exp \Bigl( \gamma \sum_{j \in \mathbb{T}_N} \Big{|} \frac{1}{2l+1} \sum_{k:|k-j|\leq l}\Bigl( \psi(\tau_k \omega)-\mathbb{E}_{\nu_{\rho(t,k/N)}}(\psi(\omega))\Bigr) \Big{|}\Bigr)
\\ & = \frac{1}{\gamma N} \log \mathbb{E}_{\tilde{\nu}_{\rho(t,.)}^N} \exp \Bigl( \gamma \sum_{i=0}^{s-1} \sum_{m=1}^{2l+d+1} \Big{|} \frac{1}{2l+1} \sum_{k:|k - (2l+d+1)i-m|\leq l}\Bigl( \psi(\tau_k \omega)-\mathbb{E}_{\nu_{\rho(t,k/N)}}(\psi(\omega))\Bigr) \Big{|}\Bigr)
\\ & \leq  \frac{1}{\gamma (2l+d+1) N} \sum_{i=0}^{s-1}\sum_{m=1}^{2l+d+1} \log \mathbb{E}_{\tilde{\nu}_{\rho(t,.)}^N} \exp \biggl( (2l+d+1)\gamma \  \times
\\ & \qquad \times \Big{|} \frac{1}{2l+1} \sum_{k:|k - (2l+d+1)i-m|\leq l}\Bigl( \psi(\tau_k \omega)-\mathbb{E}_{\nu_{\rho(t,k/N)}}(\psi(\omega))\Bigr) \Big{|}\biggr)
 \\ & \leq  \frac{1}{\gamma (2l+1) N}  \sum_{i=0}^{s-1} \sum_{m=1}^{2l+d+1} \log \mathbb{E}_{\tilde{\nu}_{\rho(t,.)}^N} \exp \biggl( 2\gamma \Big{|} \sum_{k : |k-(2l+d+1)i-m| \leq l} \Bigl( \psi(\tau_k \omega)
 -\mathbb{E}_{\nu_{\rho(t,k/N)}}(\psi(\omega))\Bigr) \Big{|}\biggr)
\\ & \text{(since} \ \frac{2l+d+1}{2l+1} \leq 2 \ \text{for sufficiently large} \ l \ \text{)}
 \\ & = \frac{1}{\gamma (2l+1) N} \sum_{j \in \mathbb{T}_N} \log \mathbb{E}_{\tilde{\nu}_{\rho(t,.)}^N} \exp \Bigl( 2\gamma  \Big{|} \sum_{k:|k-j|\leq l}\Bigl( \psi(\tau_k \omega)-\mathbb{E}_{\nu_{\rho(t,k/N)}}(\psi(\omega))\Bigr) \Big{|}\Bigr).
\end{aligned} \]
This uses the independence of integrals which are $2l+d+1$ sites apart.

\noindent
\textbf{Step 5.}
\[ \begin{aligned} & \frac{1}{\gamma (2l+1) N}  \sum_{j \in \mathbb{T}_N} \log \mathbb{E}_{\tilde{\nu}_{\rho(t,.)}^N} \exp \Bigl( 2 \gamma  \Big{|} \sum_{k:|k-j|\leq l}\Bigl( \psi(\tau_k \omega)-\mathbb{E}_{\nu_{\rho(t,k/N)}}(\psi(\omega))\Bigr) \Big{|}\Bigr)
\\ & \leq \frac{1}{\gamma(2l+1)}\int_{\mathbb{T}}\log \mathbb{E}_{\nu_{\rho(t,u)}}\exp \Bigl( 2 \gamma \Big{|} \sum_{k=0}^{2l} \Bigl( \psi(\tau_k \omega)-\mathbb{E}_{\nu_\rho(t,u)}(\psi(\omega)) \Bigr) \Big{|}\Bigr) \ \text{d}u 
\\ & \quad + \Bigg{|} \frac{1}{\gamma(2l+1)}\int_{\mathbb{T}}\log \mathbb{E}_{\nu_{\rho(t,u)}}\exp \Bigl( 2 \gamma \Big{|} \sum_{k=0}^{2l} \Bigl( \psi(\tau_k \omega)-\mathbb{E}_{\nu_\rho(t,u)}(\psi(\omega))\Bigr) \Big{|}\Bigr) \ \text{d}u 
\\ & \qquad  - \frac{1}{\gamma (2l+1) N} \sum_{j \in \mathbb{T}_N} \log \mathbb{E}_{\nu_{\rho(t,.)}^N} \exp \Bigl( 2 \gamma  \Big{|} \sum_{k:|k-j|\leq l}\Bigl( \psi(\tau_k \omega)-\mathbb{E}_{\nu_{\rho(t,k/N)}}(\psi(\omega))\Bigr) \Big{|}\Bigr) \Bigg{|} \\ & \quad
+ \Bigg{|} \frac{1}{\gamma (2l+1) N} \sum_{j \in \mathbb{T}_N} \log \mathbb{E}_{\tilde{\nu}_{\rho(t,.)}^N} \exp \Bigl( 2 \gamma  \Big{|} \sum_{k:|k-j|\leq l}\Bigl( \psi(\tau_k \omega)-\mathbb{E}_{\nu_{\rho(t,k/N)}}(\psi(\omega))\Bigr) \Big{|}\Bigr) \\ & \qquad
- \frac{1}{\gamma (2l+1) N} \sum_{x \in \mathbb{T}_N} \log \mathbb{E}_{\nu_{\rho(t,.)}^N} \exp \Bigl( 2 \gamma  \Big{|} \sum_{k:|k-j|\leq l}\Bigl( \psi(\tau_k \omega)-\mathbb{E}_{\nu_{\rho(t,k/N)}}(\psi(\omega))\Bigr) \Big{|}\Bigr) \Bigg{|}
\end{aligned} \]
The error term on lines 3 and 4 of the previous display consists of Riemann sums and converges to zero due to piecewise continuity. 
The term on lines 5 and 6 is controlled by the same coupling argument that was used in Section \ref{sc:4.4}. 

\noindent
\textbf{Step 6.} As in \cite[p.118]{landim}, since $\psi$ is bounded the inequalities $e^x \leq 1 + x + \frac{1}{2}x^2 e^{|x|}$ and $\log (1+x) \leq x$ can be used to bound
$$\frac{1}{\gamma(2l+1)}\int_{\mathbb{T}}\log \mathbb{E}_{\nu_{\rho(t,u)}}\exp \Bigl( 2 \gamma \Big{|} \sum_{k=0}^{2l} \Bigl( \psi(\tau_k \omega)-\mathbb{E}_{\nu_\rho(t,u)}(\psi(\omega)) \Bigr) \Big{|}\Bigr) \ \text{d}u $$
by \[
\begin{aligned} \frac{1}{\gamma(2l+1)}\int_{\mathbb{T}} &  \Bigg{\{} 2\gamma \mathbb{E}_{\nu_{\rho(t,u)}} \Bigl(  \Big{|} \sum_{k=0}^{2l} \Bigl( \psi(\tau_k \omega)-\mathbb{E}_{\nu_\rho(t,u)}(\psi(\omega)) \Big{|} \Bigr) \Bigr)
\\ & + 8\gamma^2 (2l+1)^2 \| \psi \|_{\infty}^2 \exp ( 4\gamma (2l+1) \| \psi \|_\infty ) \Bigg{\}} \ \text{d}u.
\end{aligned} \]
Then fix $\varepsilon > 0$ and let $\gamma = \cfrac{\varepsilon}{2(2l+1)}$, so that $\gamma$ is a function of $l$. We send $l \to \infty$ and then $\varepsilon \to 0$. The first term inside the bracket converges to zero by the law of large numbers and the second term has the form $2\varepsilon^2 \| \psi \|_\infty^2 \exp (2\varepsilon \| \psi \|_\infty)$ and converges to zero as $\varepsilon \to 0$.
\end{proof}

\section{One block estimates}\label{sc:5}

The one block estimate is proved for the general case of a non-translation-invariant graph containing a mixture of figures `1', `2' and `3', following and adapting the steps of the proof of Lemma 5.3.1 in Kipnis and Landim \cite{landim}. This one block estimate is Lemma \ref{lm:3} in Section \ref{sc:4.3}.

\begin{lemma}
 $$\limsup_{l \to \infty} \limsup_{N \to \infty} \frac{1}{T_N} \int_{\Omega_N} \sum_{j=1}^{T_N} \Big{|} g_j^l(\omega)-2\Phi(\omega_j^l) \Big{|} f_t^N \text{d}\nu^N = 0$$
 where $f_t^N = \cfrac{\text{d}\mu_t^N}{\text{d}\nu^N}$ and $g_j^l(\omega)= \frac{1}{4l+2}\sum_{|k-j| \leq l}(g(\omega_{x_{k},1})+g(\omega_{x_{k},-1}))$.
\end{lemma}

\begin{proof}
 We follow the steps of the proof of Lemma 5.3.1 in Kipnis and Landim \cite[p.82]{landim}, adapting it where necessary.

\noindent
\textbf{Step 1.} The expression is split into parts with a high and low particle density, and the part with high density is controlled. Let $K > 0$. Following the argument in \cite[pp.84-85]{landim}, the expression can be rewritten as 
\[ \begin{aligned}
\int_{\Omega_N} \sum_{j=1}^{T_N} \Big{|} g_j^l(\omega)-2\Phi(\omega_j^l) \Big{|} h^N \text{d}\nu^N & = \int_{\Omega_N} \sum_{j=1}^{T_N} \Big{|} g_j^l(\omega)-2\Phi(\omega_j^l) \Big{|}\textbf{1}\{ |\omega_j^l  | \leq K \} f_t^N \text{d}\nu^N \\ & \quad
+ \int_{\Omega_N} \sum_{j=1}^{T_N} \Big{|} g_j^l(\omega)-2\Phi(\omega_j^l) \Big{|} \textbf{1}\{ |\omega_j^l | > K \} f_t^N \text{d}\nu^N. \end{aligned} \]

To prove the lemma, it suffices to show that for all $a > 0$,
\begin{equation} \limsup_{l \to \infty} \limsup_{N \to \infty} \frac{1}{N}  \sum_{j=1}^{T_N} \int_{\Omega_N} \Bigl( \Big{|} g_j^l(\omega)-2\Phi(\omega_j^l) \Big{|} - a\omega_j^l \Bigr) f_t^N \text{d} \nu^N \leq 0, \label{eq:14} \end{equation}
 since $a$ can be made arbitrarily small and due to Lemma \ref{lm:2} the part we subtract is bounded, that is $$\frac{1}{N}\int_{\Omega_N} \sum_{j=1}^{T_N} \omega_j^l f_t^N \text{d}\nu^N \leq C.$$

Under the condition \textbf{(SLG)} on $g$, for any $b > 0$, $g(k) \leq C(b)+bk$. Taking the expectation with respect to $\nu_{k}^N$, $\Phi(k) \leq C(b) + bk$. Therefore $$\Phi(\omega_0^l) \leq C(b) + b\omega_0^l.$$ This step is valid since although $\omega_0^l$ is not necessarily an integer, the function $\Phi$ is smooth and increasing. Therefore, letting $b = \frac{1}{4}a$, 
$$|g_0^l(\omega) - 2\Phi(\omega_0^l)| \leq 2\Phi(\omega_0^l) + |g_0^l(\omega)| \leq 2C(a/4)+\frac{1}{2}a \omega_0^l.$$
 So $$|g_0^l(\omega)-2\Phi(\omega_0^l)| - a\omega_0^l \leq 0$$ if $\omega_0^l \geq \frac{4}{a}C(a/4)$. Hence the integrand in 
\eqref{eq:14} is negative when $\omega_0^l \geq \frac{4}{a}C(a/4)$, so this case can be ignored leaving the part with bounded density. We can take $K = \frac{4}{a}C(a/4)$ in the equation above.
Hence it suffices to consider the term with $\{ |\omega_j^l  | \leq K \}$ for a sufficiently large value of $K$.

\noindent
\textbf{Step 2.}  We aim to rewrite 
$$ \frac{1}{T_N} \int_{\Omega_N} \sum_{j = 1}^{T_N} \Big{|} g_j^l(\omega) - 2\Phi(\omega_j^l)\Big{|} \textbf{1}\{ |\omega_j^l |  \leq K \} f_t^N \text{d} \nu^N$$
as a sum of translated terms centred at the origin. 
The $2l+1$ tiles surrounding the pairs $(x_j, 1)$, $(x_j,-1)$ have differing numbers of vertices. For a fixed $l$, we enumerate all possible shapes a sequence of $2l+1$ tiles may have and and let $S_l^m$ denote the $m-$th shape.
We split the sum into components corresponding to the centres $(x_j,  1)$ and $(x_j,-1)$ with each surrounding tile shape. Let $C_{N,l}^m$ denote the number of centres with surrounding tile shape $S_l^m$ in the $N-$th graph, and let $A_{N,l}^m$ denote the set of sites in $\mathbb{T}_N$ with surrounding shape $S_l^m$. We define the Radon-Nikoydym derivatives
$$\bar{f}^N = \cfrac{1}{T_N}\sum_{j=1}^{T_N} f_t^N (\tau_j \omega),$$ the mean of all translated Radon-Nikodym derivatives $f_t^N(\tau_j \omega)$, and  
 $$\bar{f}_{S_l^m}^N(\omega) = \cfrac{ 1 }{C_{N,l}^m} \sum_{j \in A_{N,l}^m} f_t^N(\tau_j \omega)$$ the mean of only the translated  $f_t^N(\tau_j \omega)$ which are centred on a chain of $2l+1$ tiles of shape $S_l^m$. 
 We also define $$\bar{f}_{l,S_l^m}^N (\omega) = \int_{i \notin S_l^m} \bar{f}_{S_l^m}^N(\omega)\ \text{d}\nu^N (\omega_i).$$
 Observe that $\bar{f}_{l,S_l^m}^N \text{d}\nu^{2l+1,m}$ is the marginal of $\bar{f}_{S_l^m}^N \text{d}\nu^N$ on the $2l+1$ tiles of shape $S_l^m$, obtained by normalising by $\nu^{2l+1,m}(\omega)$. Let $\Omega_l^m$ be the state space of the graph restricted to $S_l^m$. Then, rearranging the sum,
\[ \begin{aligned} & \frac{1}{T_N} \int_{\Omega_N} \sum_{j = 1}^{T_N} \Big{|} g_j^l(\omega) - 2\Phi(\omega_j^l)\Big{|} \textbf{1}\{ |\omega_j^l |  \leq K \} f_t^N \text{d} \nu^N \\ & = \sum_{S_l^m} \cfrac{ C_{N,l}^m}{T_N} \int_{\Omega_N} \Big{|} g^{l}(\omega_0) - 2\Phi(\omega_0^{l}) \Big{|} \textbf{1}\{ |\omega_0^l |  \leq K \} \bar{f}_{S_l^m}^N \text{d}\nu^N  \\ & = \sum_{S_l^m} \cfrac{ C_{N,l}^m }{T_N} \int_{\Omega_l^m} \Big{|} g^{l}(\omega_0) - 2\Phi(\omega_0^{l}) \Big{|} \textbf{1}\{ |\omega_0^l |  \leq K \} \bar{f}_{l,S_l^m}^N \text{d}\nu^{2l+1,m} . \end{aligned} \]
This is a sum of integrals over chains of $2l+1$ tiles of varying lengths. However, note that each chain of tiles does not necessarily form a torus.

\noindent
\textbf{Step 3.}
Let $D_N(h)$ be the Dirichlet form on $\mathbb{T} \times \{-1,1\}$, defined as $$D_N(h) = \sum_{i \in \mathbb{T}_N} \Bigl( I_i^{\text{type}(i),1}(h) + I_i^{\text{type}(i),2}(h) + I_i^{\text{type}(i),3}(h)  + I_i^{\text{type}(i),4}(h) \Bigr),$$ where  $\text{type}(i) \in \{ 1, 2, 3 \}$ is the connecting figure type between $i$ and $i+1$, and
\[ \begin{aligned} &
I_i^{1,1}(h) = \frac{1}{2}\int_{\Omega_N} g(\omega_{i,1}) \Bigl( \sqrt{h(\omega^{(i,1),(i+1,1)})}-\sqrt{h(\omega)}\Bigr)^2 \text{d}\nu^N \\ &
I_i^{1,2}(h) = \frac{1}{2}\int_{\Omega_N} g(\omega_{i,-1}) \Bigl( \sqrt{h(\omega^{(i,-1),(i+1,-1)})}-\sqrt{h(\omega)}\Bigr)^2 \text{d}\nu^N  \\ &
I_i^{1,3}(h) = \frac{1}{2}\int_{\Omega_N} g(\omega_{i+1,1}) \Bigl( \sqrt{h(\omega^{(i+1,1),(i,-1)})}-\sqrt{h(\omega)}\Bigr)^2 \text{d}\nu^N \\ &
I_i^{1,4}(h) = \frac{1}{2}\int_{\Omega_N} g(\omega_{i+1,-1}) \Bigl( \sqrt{h(\omega^{(i+1,-1),(i,1)})}-\sqrt{h(\omega)}\Bigr)^2 \text{d}\nu^N  \\ &
I_i^{2,1}(h) = \frac{1}{2}\int_{\Omega_N} g(\omega_{i,1}) \Bigl( \sqrt{h(\omega^{(i,1),(i+1,1)})}-\sqrt{h(\omega)}\Bigr)^2 \text{d}\nu^N \\ &
I_i^{2,2}(h) = \frac{1}{2}\int_{\Omega_N} g(\omega_{i,-1}) \Bigl( \sqrt{h(\omega^{(i,-1),(i+1,1)})}-\sqrt{h(\omega)}\Bigr)^2 \text{d}\nu^N  \\ &
I_i^{2,3}(h) = \frac{1}{2}\int_{\Omega_N} g(\omega_{i+1,-1}) \Bigl( \sqrt{h(\omega^{(i+1,-1),(i,1)})}-\sqrt{h(\omega)}\Bigr)^2 \text{d}\nu^N \\ &
I_i^{2,4}(h) = \frac{1}{2}\int_{\Omega_N} g(\omega_{i+1,-1}) \Bigl( \sqrt{h(\omega^{(i+1,-1),(i,-1)})}-\sqrt{h(\omega)}\Bigr)^2 \text{d}\nu^N \\ &
I_i^{3,1}(h) = \frac{1}{2}\int_{\Omega_N} g(\omega_{i,1}) \Bigl( \sqrt{h(\omega^{(i,1),(i+1,1)})}-\sqrt{h(\omega)} \Bigr)^2 \text{d}\nu^N \\ &
I_i^{3,2}(h) = \frac{1}{2}\int_{\Omega_N} g(\omega_{i,1}) \Bigl( \sqrt{h(\omega^{(i,1),(i+1,-1)})}-\sqrt{h(\omega)} \Bigr)^2 \text{d}\nu^N \\ &
I_i^{3,3}(h) = \frac{1}{2}\int_{\Omega_N} g(\omega_{i+1,1}) \Bigl( \sqrt{h(\omega^{(i+1,1),(i,-1)})}-\sqrt{h(\omega)} \Bigr)^2 \text{d}\nu^N \\ &
I_i^{3,4}(h) = \frac{1}{2}\int_{\Omega_N} g(\omega_{i+1,-1}) \Bigl( \sqrt{h(\omega^{(i+1,-1),(i,-1)})}-\sqrt{h(\omega)} \Bigr)^2 \text{d}\nu^N .
\end{aligned} \]
The notation $\omega^{x,y}$ denotes the configuration $\omega$ with one less particle at $x$ and one more particle at $y$. Note that $\omega^{x,y} \neq \omega^{y,x}$, and $I_i^{j,k}$ denotes the $k$-th term of the part of the Dirichlet form between $i$ and $i+1$, where the connecting figure is of type $j$. 
Fix a shape $S_l^m$. The marginal of $D_N$ on $2l+1$ tiles of shape $S_l^m$ centred on $(0, 1)$ and $(0,-1)$ is defined by

$$D_{S_l^m}(h) = \sum_{i = x_{-l}-2}^{x_l+1} I_i^{\text{type}(i),1}(h) + I_i^{\text{type}(i),2}(h) + I_i^{\text{type}(i),3}(h)  + I_i^{\text{type}(i),4}(h),$$ where any terms at the endpoints $x_{-l}-2$ and $x_l + 1$ corresponding to vertices not included in the $2l+1$ tiles are set equal to zero.

We will now find an upper bound for $D_{S_l^m}\Bigl( \bar{f}_{l,S_l^m}^N \Bigr)$. 
We have that \[ \begin{aligned} D_{S_l^m} \Bigl( \bar{f}_{l,S_l^m}^N \Bigr) 
& \leq D(\bar{f}_{S_l^m}^N) \ \text{(since each} \ I_i^{j,k}(\bar{f}_{l,S_l^m}^N) \leq I_i^{j,k}(\bar{f}_{S_l^m}^N)   \ \text{by convexity)}
\\ & = D \Bigl( \frac{1}{C_{N,l}^m} \sum_{j \in A_{N,l}^m} f_t^N(\tau_j \omega) \Bigr)
\\ & \leq \frac{1}{C_{N,l}^m} \sum_{j \in A_{N,l}^m} D \Bigl( f_t^N (\tau_j \omega) \Bigr) \ \text{(convexity)}
\\ & \leq \frac{1}{C_{N,l}^m} \sum_{j \in \mathbb{T}_N} D \Bigl( f_t^N(\tau_j \omega) \Bigr) \ \text{(since all terms are nonnegative).}
 \end{aligned} \]
 It is not possible to simplify the final line further since $f_t^N$ is not translation invariant.

\noindent
\textbf{Step 4.} 
For a fixed shape $S_l^m$, we have that
\begin{equation*} \begin{aligned} & 
\limsup_{N \to \infty} \int_{\Omega_l^m} \Big{|} g^{l}(\omega_0) - 2\Phi(\omega_0^{l}) \Big{|} \textbf{1}\{ |\omega_j^l |  \leq K \} \bar{f}_{l,S_l^m}^N(\omega) \text{d}\nu^{2l+1,m} \\ & \quad \leq \sup_{f_{S_l^m}: D_{S_l^m}(f_{S_l^m})=0} \int_{\Omega_l^m} \Big{|} g_0^l(\omega) - 2\Phi(\omega_0^l) \Big{|} \textbf{1}\{ |\omega_0^l | \leq K \} f_{S_l^m} (\omega)\text{d}\nu^{2l+1,m} \end{aligned}
\end{equation*} by the subsequence argument given in Kipnis and Landim \cite[p.88]{landim}. This argument relies on the fact that $\lim_{N \to \infty} D_{S_l^m} \Bigl( \bar{f}_{l,S_l^m}^N \Bigr) = 0$. To establish this, we first note that by Step 3,  $D_{S_l^m} \Bigl( \bar{f}_{l,S_l^m}^N \Bigr) \leq \frac{1}{C_{N,l}^m} \sum_{j \in \mathbb{T}_N} D\Bigl( f_t^N(\tau_j \omega) \Bigr)$ and apply the bound $D_N(f_t^N(\tau_j \omega)) \leq \frac{C}{N}$ for the Radon-Nikodym derivative $f_t^N$ after integrating and averaging over $t$ \cite[p.81]{landim}. This bound holds for a constant $C$ uniformly over all translations $\tau_j$ for $j \in \mathbb{T}_N$ due to the translation invariance of $\nu^N$. Finally we use the assumption \textbf{(G2)} that $C_{N,l}^m \to \infty$ as $N \to \infty$.

Then let $p(S_l^m) = \lim_{N \to \infty} \cfrac{ C_{N,l}^m }{T_N}$ be the limiting proportion of chains of $2l+1$ tiles of shape $S_l^m$. By \textbf{(G1)} this quantity exists for each $S_l^m$ and $\sum_{m} p(S_l^m) = 1$. Then
\begin{equation*} \begin{aligned} &
\limsup_{N \to \infty} \sum_{S_l} \cfrac{ C_{N,l}^m }{T_N} \int_{\Omega_l^m} \Big{|} g_0^{l}(\omega) - 2\Phi(\omega_0^{l}) \Big{|} \textbf{1}\{ |\omega_j^l |  \leq K \} \bar{f}_{l,S_l^m}^N(\omega) \text{d}\nu^{2l+1,m} 
\\ & \quad \leq \sum_{S_l^m} p(S_l^m) \sup_{f_{S_l^m}: D_{S_l^m}(f_{S_l^m})=0} \int_{\Omega_l^m} \Big{|} g_0^l(\omega) - 2\Phi(\omega_0^l) \Big{|} \textbf{1}\{ |\omega_0^l | \leq K \} f_{S_l^m} (\omega)\text{d}\nu^{2l+1,m} \end{aligned}
\end{equation*} by weak convergence, and using that the integrand is bounded due to the $\{ |\omega_j^l  | \leq K \}$ factor.

\noindent
\textbf{Step 5.} We note that if $D_{S_l^m}(f_{S_l^m}(\omega)) = 0$ then $f_{S_l^m}$ is constant for each number of particles on the $2l+1$ tiles. Let

$$w(S_l^m,s)  = p(S_l^m) \int_{\Omega_l^m} \textbf{1}\Big{\{} \sum \omega_j = s \Big{\}} f_{S_l^m} (\omega)\text{d}\nu^{2l+1,m}.$$
Note that $\sum_{S_l^m} \sum_{s=0}^\infty w(S_l^m,s) =  \sum_{S_l^m}p(S_l^m) =  1$ since each $f_{S_l^m}$ is a Radon-Nikodym derivative. Let $\nu_s^{2l+1,m}$ be the conditional distribution of $\nu^{2l+1,m}$ given that there are a total of $s$ particles on the $2l+1$ tiles. Let $N_{S_l^m}$ denote the number of vertices in $S_l^m$. We have that
\[ \begin{aligned} &
\sum_{S_l^m} p(S_l^m)  \int_{\Omega_l^m} \Big{|} g_0^l(\omega) - 2\Phi(\omega_0^l) \Big{|} \textbf{1}\{ |\omega_0^l | \leq K \} f_{S_l^m} (\omega)\text{d}\nu^{2l+1,m} \\ & \quad
=  \sum_{S_l^m} \sum_{s=0}^{KN_{S_l^m}} w(S_l^m,s) \int_{\Omega_l^m} \Big{|} g_0^l(\omega) - 2\Phi \Bigl( \frac{s}{N_{S_l^m}} \Bigr) \Big{|} \text{d}\nu_s^{2l+1,m} 
\\ & \quad
\leq \sup_{S_l^m, s \in [0,KN_{S_l^m}]}  \int_{\Omega_l^m} \Big{|} g_0^l(\omega) - 2\Phi \Bigl( \frac{s}{N_{S_l^m}} \Bigr) \Big{|} \text{d}\nu_s^{2l+1,m} .
\end{aligned} \]

\noindent
\textbf{Step 6.} For fixed $k \in \mathbb{N}$, $S_l^m$ and $s \in \{0, \cdots, KN_{S_l^m} \}$ fixed,
\[ \begin{aligned} &
\int_{\Omega_l^m} \Big{|} g_0^l(\omega) - 2\Phi \Bigl( \frac{s}{N_{S_l^m}} \Bigr) \Big{|} \text{d}\nu_s^{2l+1,m} \\ & = \int_{\Omega_l^m} \frac{1}{2l+1} \Big{|} \sum_{j = -l}^{l} (g(\omega_{x_j,1})+g(\omega_{x_j,-1})) - 2\Phi \Bigl( \frac{s}{N_{S_l^m}} \Bigr) \Big{|} \text{d}\nu_s^{2l+1,m} \\ &
= \int_{\Omega_l^m} \frac{1}{2l+1} \Big{|} \sum_{j = -l}^{l} g_j^k(\omega) - 2\Phi \Bigl( \frac{s}{N_{S_l^m}} \Bigr)  \Big{|} \text{d}\nu_s^{2l+1,m} \\ &
\text{(where the terms of} \ g_j^k(\omega) \ \text{wrap around the ends of} \ S_l^m \ \text{ in the sense that if} \ |j|>l \ \text{then}  \\ & \quad g(\omega_{x_j,1}) \ \text{is understood as} \ g(\omega_{x_{j+(2l+1)n},1}) \  \text{ where} \ n \ \text{is such that} \ j+(2l+1)n \in [-l,l] \ \text{)}
\\ &
\leq \int_{\Omega_l^m} \frac{1}{2l+1} \sum_{j = -l}^{l} \Big{|} g_j^k(\omega) - 2\Phi \Bigl( \frac{s}{N_{S_l^m}} \Bigr) \Big{|} \text{d}\nu_s^{2l+1,m}
\\ &  =
\int_{\Omega_l^m} \Big{|} g_0^k(\omega) - 2\Phi \Bigl( \frac{s}{N_{S_l^m}} \Bigr) \Big{|} \text{d}\nu_s^{2l+1,m} \  \text{(by exchangeability).} 
\end{aligned} \]
Hence
\begin{equation*} \begin{aligned} &
\sup_{S_l^m, s \in [0,KN_{S_l^m}]}   \int_{\Omega_l^m} \Big{|} g_0^l(\omega) - 2\Phi \Bigl( \frac{s}{N_{S_l^m}} \Bigr) \Big{|} \text{d}\nu_s^{2l+1,m}  \\ & \quad \leq \sup_{S_l^m, s \in [0,KN_{S_l^m}]}  \int_{\Omega_l^m} \Big{|} g_0^k(\omega) - 2\Phi \Bigl( \frac{s}{N_{S_l^m}} \Bigr) \Big{|} \text{d}\nu_s^{2l+1,m}. \end{aligned} \end{equation*}
Let $u \in \mathbb{T}$, with the corresponding expected particle density $\rho(t,u)$ at time $t$. 
Consider a sequence of particle numbers for each torus length $l$ and corresponding set of shapes $S_l^m$ such that the particle density converges to $\rho(t,u)$. That is, $$ \lim_{l \to \infty} \sup_{m} \Big{|} \frac{s}{N_{S_l^m}} - \rho(t,u) \Big{|} = 0.$$
Applying the equivalence of ensembles, the expected value of $|g_0^l(\omega)-2\Phi(s/N_{S_l^m})|$ under $\nu_s^{2l+1,m}$, the marginal conditional on there being a total of $s$ particles occupying the $N_{S_l^m}$ sites, converges to its expected value under $\nu_{\rho(t,u)}^k$
and
\[ \begin{aligned} & \limsup_{l \to \infty} \ \sup_{S_l^m} \int_{\Omega_l^m} \Big{|} g_0^k(\omega) - 2\Phi \Bigl( \frac{s}{N_{S_l^m}} \Bigr) \Big{|} \text{d}\nu_s^{2l+1,m} \\ & \leq \sup_{u: \rho(t,u) \in [0,K]} \int_{\Omega_{2k+1}} \Big{|} g_0^k(\omega) - 2\Phi (\rho(t,u)) \Big{|} \text{d}\nu_{\rho(t,u)}^k \\ & \quad \to 0 \ \text{as} \ k \to \infty \end{aligned} \]
as required.
\end{proof}

\section*{Acknowledgements}
The authors would like to thank B\'alint T\'oth for his suggestion to study the hydrodynamic limit of a totally asymmetric interacting particle system in a random environment using the relative entropy method, Stefano Olla for his suggestion to consider a zero range process, and anonymous Referees for valuable suggestions to improve the manuscript. Felix Maxey-Hawkins was supported by an EPSRC studentship and M\'arton Bal\'azs was partially supported by the EPSRC EP/R021449/1 Standard Grant of the UK. This study did not involve any underlying data.

\end{document}